\documentclass{ps}

\usepackage[colorlinks = true,
            linkcolor = blue,
            urlcolor = amber,
            citecolor = blue,
            anchorcolor = blue]{hyperref}

\usepackage[utf8]{inputenc}
\usepackage[T1]{fontenc}

\newcommand{\e}[0]{\mathbf{e}}

\newcommand{\iii}[0]{\mathbf{i}}
\newcommand{\rr}[0]{\mathbf{r}}
\newcommand{\kkappa}[0]{\boldsymbol{\kappa}}

\newcommand{\PP}[0]{\mathbf{P}}

\newcommand{\s}[0]{\mathbf{s}}
\newcommand{\C}[0]{\mathbf{C}}
\newcommand{\X}[0]{\mathbf{X}}

\newcommand{\trev}[0]{\overleftarrow}

\newcommand{\EE}{\mathbb E}

\newcommand{\diag}{\mathbf{diag}}

\usepackage{graphicx}
\usepackage{soul} 
\usepackage{tikz}

\usepackage{amssymb}
\usepackage{MnSymbol}

\newtheorem{deff}{Definition}[section]

 \newtheorem{remark}[deff]{Remark}

\begin{document}
\title{Antiduality and M\"obius monotonicity: Generalized Coupon Collector Problem}\thanks{Work supported by NCN Research Grant DEC-2013/10/E/ST1/00359}
\runningtitle{Antiduality and M\"obius monotonicity}
%
\author{Paweł Lorek}\address{Mathematical Institute, University of Wrocław, pl. Grunwaldzki 2/4, 50-384 Wrocław, Poland. Email: Pawel.Lorek@math.uni.wroc.pl}
%

%
\begin{abstract}
For a given absorbing  Markov chain $X^*$ on a finite state space, a chain $X$ is a 
sharp antidual of $X^*$ if the fastest strong stationary time of $X$ is equal, in distribution,
 to the  absorption time  of $X^*$. In this paper we show a systematic way of finding 
 such an antidual based on some partial ordering of the state space.
 We use a theory of strong stationary duality developed recently for M\"obius monotone Markov chains. 
 We give several sharp antidual chains  for Markov chain corresponding to a generalized coupon collector problem. 
 As a consequence - utilizing known results on a limiting distribution of the 
 absorption time - we indicate a separation cutoff (with its window size) in several 
 chains.
 We also present a chain which (under some conditions) has a prescribed stationary distribution 
 and its fastest strong stationary time is distributed as a prescribed mixture of sums of geometric random variables.
 \end{abstract}
%
%
\maketitle

\section{Introduction}
Strong stationary times (SST) are a probabilistic tool for bounding a rate of convergence to stationarity for Markov chains.
Aldous and Diaconis \cite{Aldous1986}, \cite{Aldous1987} gave several  examples of chains where SST were found \textsl{ad hoc}. 
Later in \cite{Diaconis1990a}   authors introduced a more systematic way 
of finding SSTs. For a given general ergodic chain 
they showed that one can  construct a so-called \textsl{strong stationary dual} (SSD)
chain,
 a chain whose  absorption time is equal to some SST,
 not only in distribution, via the coupling of the chain 
 with its SSD which is presented in  \cite{Diaconis1990a}.
 Moreover, they proved that there 
always exists \textsl{sharp} SSD, in the sense that its corresponding SST is stochastically the smallest,
in which  case it is called the \textsl{fastest strong stationary time} (FSST).\par
Their construction for general chains is purely theoretical (it involves the knowledge of the distribution 
of the chain at each step). However, they give a detailed recipe on how to construct such SSD assuming that 
the time reversed chain is stochastically monotone w.r.t. linear ordering. 
In particular, they consider birth and death chain, for which SST has the same distribution as absorption time in a dual chain,
which turns out to be an absorbing birth and death chain. 
They also show that assuming that time reversed chain is stochastically monotone one can always 
construct set-valued SSD  (see their Section 3.4 ``greedy construction of a set-valued dual'').
In this paper we actually start with some absorbing chain and show that it is a sharp SSD of 
a class (which we indicate) of ergodic chains. We exploit the results from  \cite{Lorek2012d}, where the
authors provided the recipe for constructing SSD 
on the same state space for chains, whose 
time reversal is M\"obius monotone w.r.t to some partial ordering of the state space. This significantly enlarges 
the class of chains for which SSD can be found. In many chains there is usually some natural underlying ordering of the state space
which is only partial. Moreover, the method yields the sharp SSD which is crucial for our applications.

\smallskip\par

Studying the rate of convergence of a chain to its stationary distribution,
one is often interested in a so-called 
\textsl{mixing time} (i.e., the time until the chain is ``close'' to its stationary distribution).
However, sometimes we can say much more than just a mixing time by showing that a so-called \textbf{cutoff} phenomenon occurs.
Roughly speaking, this phenomenon  describes a sharp transition in the convergence of the chain to its stationary distribution over a 
negligible period of time (cutoff window). There are two most commonly studied phenomena: separation cutoff and total variation cutoff,
which differ in a distance used to measure the convergence (separation vs. total variation distance).

\par 

\medskip\par 

The total variation cutoff was first shown for a \textsl{random transposition} card shuffling in \cite{Diaconis1981a}.
The name comes from \cite{Aldous1986}, where the authors showed that  a
\textsl{top-to-random} card shuffling exhibits a total variation cutoff.
A separation cutoff has recently been studied in few contexts. For example: in 
\cite{Diaconis2006} authors gave if and only if conditions for
the existence of  a
separation cutoff for birth and death chains (they 
use duality theory to convert convergence 
rates to hitting times and Keilson's representation of 
first hitting times) -- they show that there is a cutoff if 
and only if the product of a spectral gap and a mixing time tends 
to infinity; this was somehow extended -- in \cite{Chen2008} 
authors show that there is a cutoff measured in $L^p$-norm ($1<p\leq \infty$) 
if and only if the the spectral gap and max-$L^p$ mixing time tends to infinity;
computation of cutoff time and window size in a variety of birth and death chains 
is given in \cite{Chen2015};
a separation cutoff for skip-free chains was given in \cite{MaoZhang2016}; some other specific chains were considered in \cite{Connor2010};
in \cite{Fill1996} author gives a formula for the separation   for \textsl{Tsetlin library} chain specifying weights for which there is and there is no separation cutoff.
Several examples of both, separation and total variation cutoffs are given in \cite{LevPerWil_mixing_sec_ed},
some characterization of total variation cutoff for lazy (i.e., with probability of staying $\geq 1/2$) chains was recently given in \cite{Basu2017}. 
In \cite{Choi2016} authors give sufficient condition for skip-free 
chains to have real eigenvalues, they use Siegmund duality -- actually
antiduality -- a type of transitions of their (anti)dual resembles 
some chains we obtain for a coupon collector problem.  It is worth mentioning that 
although a sequence of birth and death chains exhibits total variation 
cutoff if and only if it exhibits separation cutoff \cite{Diaconis2006},
\cite{Ding2010}, it is not the case (in general) for other chains, 
as shown in \cite{Hermon2016}.

%
%
%
%
%
%
%
%
%
%
%
%
%
%
%
%
%
%
%

\medskip\par 

As mentioned before, FSST is equal in distribution to the absorption time of the sharp SSD chain. 
Thus, there is a close relation between a sharp SSD and a separation cutoff.  Roughly speaking, this cutoff can 
be studied by studying the limiting distribution of the absorption time of the SSD. 
This can be extremely difficult task. However, since  examples of chains with proven separation cutoff are always welcome,
we can reverse the procedure: starting with some already absorbing chain we can try to find an ergodic \textbf{sharp antidual} chain (or even a class of such antidual chains).
Such an approach was considered in \cite{Fill2009} in a context of birth and death chains only.
A connection between a separation cutoff and a coupon collector problem (including some generalizations,
e.g., sampling $k>1$ different coupons at a time) was given in  \cite{Pak2001}.

Using this approach we will indicate a separation cutoff time and a window size in several examples of chains utilizing (nontrivial) results for 
the limiting distribution of the absorption time in some generalizations of the classical coupon collector problem. 
That is why we need a recipe for sharp antidual chains, what will be given based on results from \cite{Lorek2012d}.
Most of the examples that follow deal with some product-type chains. It is however worth noting that 
taking a product of chains where each chain exhibits a cutoff does not have to yield a chain (on a product space)
exhibiting a cutoff. Such an example was recently given in \cite{Lacoin2016}.

\medskip\par 
The absorption time of many absorbing chains is distributed as a mixture of sums of geometric random variables with parameters 
being the eigenvalues of the transition matrix. E.g, the absorption time 
of discrete time birth and death chain starting at the minimal state with the
maximal one being absorbing is  distributed as a sum of geometric random variables with such parameters, provided the chain is
stochastically monotone. The result is usually attributed to Karlin and McGregor \cite{Karlin1959} or Keilson \cite{Keilson71b}.
Fill \cite{Fill2009} gave a   stochastic proof of this result using also the theory of SSD (the result was simultaneously obtained in 
\cite{Diaconis2009b}),
later it was extended to skip-free Markov chains in Fill \cite{Fill2009a}.
Miclo \cite{Miclo2010a} showed that for large class of  absorbing chains on a finite state space,
the absorption time is distributed as a mixture of sums of geometric random variables.
A natural question arises: \textsl{Given a mixture of sums of geometric random variables and some distribution $\pi$ can 
we find an ergodic chain whose stationary distribution is $\pi$ and whose FSST is equal in distribution to this mixture}?
Or, a special case of the question, \textsl{Given some distribution $\pi$ can be construct an ergodic chain whose 
stationary distribution is $\pi$ having deterministic FSST}? We provide positive answers to both questions
(some assumptions on distributions are needed). In particular, we present two ergodic chains on completely different 
state spaces having the same FSST.

\medskip\par 
The main goals of the paper are: i) we give a systematic way (based on partial ordering of the state space and M\"obius monotonicity)
for finding a class of sharp antidual chains; ii) we present nontrivial antidual chains related to some generalizations of coupon collector problem and,
as a consequence, we show cutoff phenomena in some cases; iii) we present a construction of a chain with prescribed FSST and prescribed stationary distribution.

 \medskip\par

There is yet another potential application which served as a motivation for the paper (however, not exploited here): Given a probability distribution $\pi$ on $\EE$,
how to simulate a sample from this distribution? Markov Chain Monte Carlo methods come with the answer: construct a chain with stationary
distribution $\pi$ and run it 
  \textsl{long enough}. The most common algorithms for such constructions are Metropolis-Hastings algorithm and Gibbs sampler
  (for studies on rate of convergence for Metropolis-Hasting algorithms see, e.g., \cite{DiaSal_metropolis98}, the cutoff 
  for Gibbs sampler for Ising model on the lattice was studied on \cite{Lubetzky2013}).
This paper suggests an alternative approach: given $\pi$ on $\EE$ find some absorbing chain on $\EE$ and then calculate sharp antidual chain having this $\pi$ as
stationary distribution. Knowing, e.g., expectation and variance of absorption time, one can quite precisely determine the number of steps needed for simulation.
Moreover,  having a sharp SSD actually can allow for a \textsl{perfect simulation} from distribution $\pi$. 
One can construct an appropriate coupling of 
the absorbing chain and its antidual, so that   stopping antidual chain when its SSD is absorbed yields an unbiased 
sample from $\pi$. The reader is referred for details to  \cite{Diaconis1990a} (Section 2.4),
\cite{FilLyz14} (Section 1.1) or \cite{LorMar_mon_req} (Section 2.3, Algorithm 4).  We want to emphasize that utilizing this was not the purpose of this paper, and the stationary distributions which appear
in most of  the examples are
of product form, which means we can easily  simulate them coordinate by coordinate.



 \medskip\par

The paper is organized as follows. 
In Section \ref{sec:prelim} we introduce preliminaries on strong stationary duality and separation cutoff.
In Section \ref{sec:mob} we recall a notion of M\"obius monotonicity and give a matrix-form proof
of the result from \cite{Lorek2012d}. In Section \ref{sec:main} we present our main results. Firstly, 
in Section \ref{sec:main_gen_mc} in Theorem \ref{twr:main} we give  a systematic way  for finding a class of sharp antidual chains.
Secondly, in Section \ref{sec:main_antidual_chains} we introduce in details the chain corresponding to the generalized coupon collector problem 
and present sharp antidual chains in Theorems 
\ref{twr:ver1} and \ref{twr:coupon_antidual_ver2}. Then, in Section \ref{sec:cutoff_results}, we proceed with presenting separation cutoff results
for some cases. In Section \ref{sec:constr_FSST} we present our results concerning construction of ergodic chain with 
prescribed stationary distribution and with prescribed FSST.
Section \ref{sec:proofs} includes main proofs. Section \ref{sec:gen_coup_col_proofs} contains  proofs of Theorems \ref{twr:ver1} and \ref{twr:coupon_antidual_ver2},
whereas Section \ref{sec:proof_FSST} contains the proof of Theorem \ref{thm:FSST_gen}.

\section{Preliminaries}\label{sec:prelim}
\subsection{Strong stationary duality}
Consider an  ergodic (i.e., irreducible and aperiodic) Markov chain 
$X\sim(\nu,\PP)$ on a finite state space 
$\EE=\{\e_1,\ldots,\e_M\}$ with an initial distribution $\nu$, 
a stationary distribution $\pi$ and
 a transition matrix $\PP$. 
 Let $\EE^*=\{\e_1^*,\ldots,\e_N^*\}$ be a state space of an absorbing Markov chain $X^*\sim(\nu^*,\PP^*)$, whose unique absorbing state 
 and unique irreducible class is denoted by $\e_N^*$. 
Define $\Lambda$, a matrix of size $N\times M$, to be a \textsl{link} if it is a stochastic matrix with the property: $\Lambda(\e_N^*,\e)=\pi(\e)$ for all $\e\in\EE$. We say that $X^*$ is 
a \textsl{strong stationary dual} (SSD) of $X$ with link $\Lambda$ if
\begin{equation}\label{eq:duality}
\nu=\nu^*\Lambda \quad \textrm{and} \quad \Lambda\PP=\PP^*\Lambda.
\end{equation}

Diaconis and Fill \cite{Diaconis1990a} prove that then the absorption time $T^*$ of $X^*$ is a so-called \textsl{strong stationary time} (SST) for $X$.
This is such a random variable $T$ that $X_T$ has distribution $\pi$ and $T$ is independent from $X_T$. 
The main application is in studying the rate of convergence of an ergodic chain to its stationary distribution, since for 
such  a random variable we always have:
$d_{TV}(\nu\PP^k,\pi)\leq sep(\nu\PP^k,\pi):=\max_{\e\in\EE}\left(1-{\nu\PP^k(\e)/\pi(\e)}\right)\leq P(T>k)$,
where $d_{TV}$ stands for \textsl{total variation distance}, and $sep$ stands for \textsl{separation}.
Note that $sep$ is not symmetric and thus is not 
a distance between probability measures.
The corresponding $T^*$  is \textbf{sharp} if $sep(\nu \PP^k,\pi) = P(T^*>k)$. In 
such a case, $T^*$ is called \textsl{the fastest strong stationary time}   for $X$,
which we denote by $T_{FSST}$.
For more details on this duality consult  
\cite{Diaconis1990a}. 
Moreover, duality relation (\ref{eq:duality}) allows  for stochastic  constructions,
see, e.g., \cite{Fill2009}, where stochastic proof for passage time distribution for birth and death chain was given.

Note that once we fix  $\EE^*$ and a link $\Lambda$, and if there exists a right-inverse of $\Lambda$, i.e., $\Lambda^{-1}$ we can simply calculate
from (\ref{eq:duality}):
$$\PP^*= \Lambda\PP\Lambda^{-1} \ \mathrm{and} \ \nu^*=  \nu \Lambda^{-1}.$$
If the resulting $\PP^*$ is a stochastic, irreducible and aperiodic matrix and $\nu^*$ is a probability distribution,  then  (it will always correspond to an absorbing  chain) we have found an SSD.
However, we can start with some already absorbing chain $\PP^*$, then find some $\EE$ and \textsl{some} probability distribution $\pi$ on $\EE$, and a link   $\Lambda$, so that 
$$\PP=\Lambda^{-1}\PP^*\Lambda \ \mathrm{and} \ \nu = \nu^*\Lambda.$$
If the resulting $\PP$ is a stochastic matrix, then $X\sim(\nu,\PP)$ is an ergodic chain with stationary distribution $\pi$, and  $T^*$ (time to absorption of $X^*$) is an
SST for $X$. 
In such a case, $X$ is called \textbf{antidual} of $X^*$.
Moreover, if we somehow know, that for some class of links  relation (\ref{eq:duality}) implies that $T^*$ is sharp (see Corollary \ref{cor:sharp}), then we can possibly find many different antiduals,
which all have the same fastest strong stationary time $T^*$, which has a phase-type distribution.  
In such a case $X$ is called a \textbf{sharp antidual} of $X^*$.
\subsection{Separation cutoff}
The forthcoming Theorem \ref{twr:main} indeed gives a 
recipe on how to construct a sharp antidual chain $X$ with  a specified stationary distribution $\pi$ given absorbing chain $X^*$,
both on the same state space. It means, that we have
\begin{equation}\label{eq:fsst}
sep(\nu\PP^k,\pi)=P(T_{FSST}>k)=P(T^*>k). 
\end{equation}
Thus, studying the distribution of $T_{FSST}$ is equivalent to study the distribution of $T^*$.
Furthermore, a separation cutoff can be studied by studying the properties of $T^*$.
In what follows, we introduce the notion of separation cutoff.
Since the definition of the cutoff involves  increasing state space, we add a subscript {\small {($d$)}} to transition matrices, distributions,
state space
and absorption time.
Suppose we have a sequence of ergodic Markov chains $X_{(d)}\sim(\nu_{(d)},\PP_{(d)})$ indexed by $d=1,2,\ldots$
Denote by $\pi_{(d)}$ the stationary distribution of $X_{(d)}$.
We say that this sequence exhibits a \textbf{separation cutoff at time } $t_d$ with a \textbf{window size} $w_d=o(t_d)$ if 
\begin{equation*}
 \begin{array}{lll}
  \displaystyle \lim_{c\to\infty}& \displaystyle \limsup_{d\to\infty} & \displaystyle sep(\nu_{(d)}\PP^{t_d+c w_d}_{(d)},\pi_{(d)})=0,\\[6pt] 
  \displaystyle \lim_{c\to\infty}& \displaystyle \liminf_{d\to\infty} & \displaystyle sep(\nu_{(d)}\PP^{t_d-c w_d}_{(d)},\pi_{(d)})=1.
 \end{array}
\end{equation*}
If the convergence to stationarity is measured in a 
total variation distance, we say about a \textbf{total variation cutoff}.
%
%
%

%
%

 \section{M\"obius monotonicity and duality}\label{sec:mob}
 In general, there is no recipe on  how to find an SSD, i.e., a triplet 
 $\EE^*, \PP^*, \Lambda$. In \cite{Diaconis1990a} authors give a recipe for a dual on the same state space $\EE^*=\EE$
 provided  that a  time reversed chain $\trev{X}$ is stochastically monotone with respect to total ordering. 
 In \cite{Lorek2012d} we give an extension of this result to state spaces 
 which are only partially ordered by $\preceq$.
 Then, provided that the time reversed chain $\trev{X}$ is \textsl{M\"obius} monotone (plus some 
 conditions on the initial distribution),
 we give a formula for a sharp SSD on the same state space $\EE^*=\EE$. 
 \medskip\par  
The M\"obius monotonicity seems to be a natural one for extension of main result from \cite{Diaconis1990a} to partially ordered state spaces.
In \cite{Lorek2016_Siegmund_duality}  we show that it is equivalent to the the 
existence of a Siegmund dual of a chain with given partial ordering.
For a linearly ordered state space, stochastic monotonicity of a chain is required for the 
existence of a Siegmund dual (see  \cite{Siegmund1976}),
and stochastic monotonicity of a time reversal is required for the existence of an SSD with a link being a truncated stationary distribution (see \cite{Diaconis1990a}).
Both results fail for non-linear orderings, since both require M\"obius monotonicity, which, in general, is different than the stochastic one.
The monotonicities  are equivalent for linear ordering.  
For more relations between these (and not only) monotonicities consult \cite{LorMar_mon_req}, and for 
applications of a Siegmund duality to some generalizations of a gambler's ruin problem consult \cite{2015Lorek_gambler}.
We will introduce this monotonicity by trying to solve (\ref{eq:duality}) with some given link $\Lambda$.
\par

\smallskip
We consider a finite state space $\EE=\{\e_1,\ldots,\e_M\}$ partially ordered by $\preceq$ such that $\e_M$ is the unique maximal state.
For a function $f:\EE\to \mathbf{R}$, by lower-case bold symbol    $\boldsymbol{f}$ we denote the row vector  $\boldsymbol{f}=(f(\e_1),\ldots,f(\e_M))$.\smallskip\par 
The idea is to find an SSD $X^*$ with a transition matrix $\PP^*$ on 
the same state space $\EE^*=\EE$ with a link, whose row corresponding to $\e$ is a stationary distribution of $X$ truncated to $\{\e\}^\downarrow:=\{\e': \e'\preceq \e\}$, i.e.,
\begin{equation}\label{eq:link0}
\Lambda(\e_i,\e_j)= {\pi(\e_j)\over \sum_{\e':\e'\preceq \e_i} \pi(\e')}\mathbf{1}(\e_j\preceq \e_i).
\end{equation}
Note that for all $\e\in\EE$ we have $\Lambda(\e_M,\e)=\pi(\e)$, as required.
For a given ordering let $\C(\e_i,\e_j)=\mathbf{1}(\e_i \preceq\e_j)$. 
 For the partial ordering we require only that the state which is absorbing for $X^*$, denoted throughout the paper by $\e_M$, is the unique maximal one 
 (i.e., $\C(\e_M,\e_j)=\mathbf{1}(\e_j=\e_M)$ for all $j$ and there is no $\e_{M_2}\neq \e_M$ such that $\C(\e_{M_2},\e_j)=\mathbf{1}(\e_j=\e_{M_2})$ for all $j$).
We always identify ordering $\preceq$ with the matrix $\C$, keeping in mind, that enumeration of states in $\C$ and $\PP$ must be the same.
Then the link can be written in a matrix form:
\begin{equation}\label{eq:link}
\Lambda = (\diag(\boldsymbol{\pi}\C))^{-1} \C^T \diag(\boldsymbol{\pi}),
\end{equation}
where $\diag(\boldsymbol{g})$ is a diagonal matrix with entries $g(\e_1),\ldots,g(\e_M)$.
The states can always be rearranged in such a way that $\C(\e_i,\e_j)=1$ implies $i\leq j$, what means that $\C$, and thus $\Lambda$, are invertible.
Often,  $\mu\equiv\C^{-1}$ is called the \textsl{M\"obius function} or 
the \textsl{M\"obius matrix}  of the partial order $\preceq$.
Solving (\ref{eq:duality}) for $\PP^*$ yields (recall	 that the transitions of time reversed chains are given by $\trev{\PP}=(\diag(\boldsymbol{\pi}))^{-1} \PP^T (\diag(\boldsymbol{\pi}))$)
\begin{eqnarray*}
\PP^*&=&\Lambda\PP\Lambda^{-1}=(\diag(\boldsymbol{\pi}\C))^{-1} \C^T \diag(\boldsymbol{\pi}) \PP (\diag(\boldsymbol{\pi}))^{-1} (\C^{T})^{-1} (\diag(\boldsymbol{\pi}\C))\\
&=&(\diag(\boldsymbol{\pi}\C) (\C^{-1} \trev{\PP} \C) (\diag(\boldsymbol{\pi}\C))^{-1})^T,
\end{eqnarray*}
which is a stochastic matrix if and only if 
each entry of $\C^{-1} \trev{\PP} \C$  is non-negative, in other words we say  that $\trev{\PP}$ is M\"obius monotone. 
This way we proved the main part of Theorem 2 of \cite{Lorek2012d}. We include it here, 
since this is a little bit different (matrix-form) proof. 
We  will restate the theorem  for completeness, introducing formal definitions of monotonicities first.
For given partial ordering $\preceq$ and any matrix $\PP$ (not necessarily stochastic) we define $\PP(\e,\{\e_j\}^\downarrow)=\sum_{\e':\e'\preceq \e_j} \PP(\e,\e')$
and  similarly $\PP(\e,\{\e_j\}^\uparrow)=\sum_{\e':\e'\succeq \e_j} \PP(\e,\e')$.
\begin{dfntn}\label{def:mob_mon}
Markov chain $X$ is M\"obius monotone if 
$\displaystyle \C^{-1}\PP \C\geq 0 $ (each entry   non-negative). In terms of transition probabilities, it means that 
\[\forall(\e_i, \e_j\in \EE) \quad \sum_{\e\succeq \e_i} \mu(\e_i,\e)\PP(\e,\{\e_j\}^\downarrow)\geq 0.\]
\end{dfntn}

\noindent 
Recall that for a M\"obius function we always have $\mu(\e_i,\e)=0$ whenever $\e_i\npreceq \e$.
\begin{dfntn}\label{def:fun_mob_mon}
A function $f:\EE\to \mathbf{R}$ is M\"obius monotone if $\mathbf{f}(\C^T)^{-1}\geq 0$ (each
entry non-negative). It means that
\[ \forall(\e_i\in\EE)\quad  \sum_{\e:\e\succeq \e_i} \mu(\e_i,\e) f(\e) \geq 0. \]
\end{dfntn}

\begin{rmrk} \rm
In Lorek, Szekli \cite{Lorek2012d} this M\"obius monotonicity (of both, function and chain) was called  \textsl{${}^\downarrow$-M\"obius monotonicity} 
(see Definitions 2.1 and 2.2 therein). 
\end{rmrk}

\begin{dfntn}\label{def:mon_up} 
$X$ is ${}^\uparrow$-M\"obius monotone if $\displaystyle (\C^T)^{-1}\PP\C^T\geq 0 $ (each entry   non-negative). 
\end{dfntn}

%
%

\begin{thrm}[Theorem 2 of \cite{Lorek2012d}]\label{twr:main_twr2_LS}
 Let $X\sim(\nu,\PP)$ be an ergodic Markov chain on a finite state space $\EE=\{\e_1,\ldots,\e_M\}$, 
partially ordered by $\preceq$, with a unique maximal state $\e_M$, and with a stationary distribution $\pi$.
Assume that 
\begin{itemize}
 \item[(i)] $g(\e)={\nu(\e)\over \pi(\e)}$ is M\"obius monotone,
 \item[(ii)] time reversed chain $\trev{X}$ is  M\"obius monotone.
\end{itemize}
Then there exists a strong stationary dual chain $X^*\sim(\nu^*,\PP^*)$ on $\EE^*=\EE$ with the following link
\begin{equation}\label{eq:Link}
 \Lambda = (\diag(\pi\C))^{-1} \C^T \diag(\pi). 
\end{equation}

\noindent 
Let $H(\e)=\sum_{\e'\preceq \e} \pi(\e')$.  The SSD chain is uniquely determined by
$$
\begin{array}{llllrlllllllll}
 \nu^*&= &   \nu\Lambda^{-1} &\textrm{ i.e.,} &  \nu^*(\e_i) & = & \displaystyle H(\e_i)   \sum_{\e:\e\succeq \e_i}\mu(\e_i,\e) g(\e),
\\[12pt]
\PP^* & = & \Lambda\PP\Lambda^{-1}, &\textrm{ i.e.,}  & \PP^*(\e_i,\e_j) & = &\displaystyle  {H(\e_j)\over H(\e_i)} \sum_{\e:\e\succeq \e_j} \mu(\e_j,\e)\trev{\PP}(\e,\{\e_i\}^\downarrow).
& \\
\end{array}
$$
\end{thrm}
\noindent
The following Corollary will play a crucial role:
\begin{crllr}\label{cor:sharp}
 The SSD constructed in Theorem \ref{twr:main_twr2_LS} is \textbf{sharp}.
\end{crllr}
\begin{proof}
The link given in  (\ref{eq:Link}) is lower-triangular, thus, by Remark 2.39 in \cite{Diaconis1990a}, the resulting SSD is sharp.
\end{proof}

\section{Main results}\label{sec:main}
\subsection{General procedure for sharp anti-dual  chains}\label{sec:main_gen_mc}
The main contribution is a systematic way of finding a sharp antidual (on the same state space $\EE=\EE^*$) chain of some given  already absorbing chain
$X^*\sim(\nu^*,\PP^*)$ with the unique absorbing state $\e_M$.
The idea is clear from the previous section: introduce some partial ordering and some distribution $\pi$ on $\EE$. Then solve $\Lambda\PP=\PP^*\Lambda$ for $\PP$ with the link given 
in (\ref{eq:Link}). If the resulting matrix is non-negative, it will be a stochastic matrix of an ergodic Markov chain $X$ with the stationary distribution $\pi$.
Moreover, changing $\pi$ and/or ordering usually will yield  a different sharp antidual. It means we can have a class of chains, 
all having the same fastest strong stationary time $T_{FFST}$.
\par
Fix some partial ordering $\preceq$ on  $\EE^*$ (expressed by $\C$) having the  unique maximal state $\e_M$ and some  distribution $\pi$ on $\EE$. For given $\PP^*$ define
$$\widehat{\PP^*} = \diag(\boldsymbol{\pi}\C) \PP^*  (\diag(\boldsymbol{\pi}\C))^{-1}.$$
With slight abuse of notation we will assume that $\widehat{\PP^*}$ is ${}^\uparrow$-M\"obius monotone meaning that $ (\C^T)^{-1}\widehat{\PP^*}\C^T\geq 0$.
Definition \ref{def:mon_up}  was stated for a Markov chain $X$ with a transition matrix $\PP$, note however that $\widehat{\PP^*}$ does not have to 
be a stochastic matrix.


\begin{thrm}\label{twr:main}
 Let $\X^* \sim (\nu^*,\PP^*)$ be an absorbing Markov chain on $\EE^*=\{\e_1,\ldots,\e_M\}$ with
 the 
unique absorbing state $\e_M$. 
Let $\mathcal{C}$ be the  class  of all partial orderings on $\EE^*$ with   $\e_M$ being unique maximal state.
Consider the class of pairs of distributions and partial orderings such that $\widehat{\PP^*}$ is ${}^\uparrow$-M\"obius monotone:
$$\mathcal{P}(\PP^*)=\left\{(\pi,\C): \C\in \mathcal{C}, \widehat{\PP^*} \textrm{\ is } {}^\uparrow\textrm{-M\"obius monotone} \right\}.$$
Then for any $(\pi,\C)\in\mathcal{P}(\PP^*)$ the chain $X\sim(\nu,\PP)$ with the  link  $\Lambda$ defined in (\ref{eq:Link}) and with
$$\nu=\nu^*\Lambda,\quad \PP=(\diag(\boldsymbol{\pi}))^{-1} (\C^T)^{-1}\widehat{\PP^*}\C^T \diag(\boldsymbol{\pi})$$
is a sharp antidual for $\PP^*$, i.e., $\PP^*$ is a sharp SSD for $\PP$.
Equivalently, $\PP=\Lambda^{-1}\PP^*\Lambda$, where, for given $\pi$ and $\C$, the link is defined in (\ref{eq:Link}).
\end{thrm}
\begin{proof}
Since $\nu^*$ is a distribution on $\EE$ and $\Lambda$ is a stochastic matrix, $\nu$ is a distribution on $\EE$. By assumption 
that $\widehat{\PP^*}$ is ${}^\uparrow$-M\"obius monotone,
the matrix $\PP$ is non-negative. We will show that $\pi$ is its stationary distribution. 
Let $\boldsymbol{\eta}=(0,\ldots,0,1)$. Last row of $\Lambda$ is equal to $\boldsymbol{\pi}$ what can be expressed as $\boldsymbol{\eta}\Lambda=\boldsymbol{\pi}$, thus 
$\boldsymbol{\eta}=\boldsymbol{\pi}\Lambda^{-1}$. 
We have
$$\boldsymbol{\pi}\PP=\boldsymbol{\pi}\Lambda^{-1}\PP^*\Lambda=\boldsymbol{\eta} \PP^*\Lambda =\boldsymbol{\eta} \Lambda = \boldsymbol{\pi}.$$

\noindent 
Now we will show that the rows of $\PP$ sum up to 1, i.e., that  $\PP (1,\ldots,1)^T=(1,\ldots,1)^T$. We have
\begin{eqnarray}\nonumber
\PP(1,\ldots,1)^T &= &(\diag \boldsymbol{\pi})^{-1} (\C^T)^{-1} \widehat{\PP^*}\C^T \diag(\boldsymbol{\pi})(1,\ldots,1)^T  \\[6pt]\nonumber
    &= &(\diag \boldsymbol{\pi})^{-1} (\C^T)^{-1} \widehat{\PP^*}\C^T \boldsymbol{\pi}^T =(\diag \boldsymbol{\pi})^{-1} (\C^T)^{-1} \widehat{\PP^*}(\boldsymbol{\pi}\C)^T  \\[6pt]\nonumber
        &= &(\diag \boldsymbol{\pi})^{-1} (\C^T)^{-1} \diag(\boldsymbol{\pi}\C) \PP^*  (\diag(\boldsymbol{\pi}\C))^{-1}(\boldsymbol{\pi}\C)^T  \\[6pt]\nonumber
        &= &(\diag \boldsymbol{\pi})^{-1} (\C^T)^{-1} \diag(\boldsymbol{\pi}\C) \PP^* (1,\ldots,1)^T   \\[6pt]\label{eq:sumup1}                
        &= &(\diag \boldsymbol{\pi})^{-1} (\C^T)^{-1} \diag(\boldsymbol{\pi}\C) (1,\ldots,1)^T \stackrel{(*)}{=} (1,\ldots,1)^T.
\end{eqnarray}
To show $(*)$ we need to show that $\sum_{\e'\in\EE^*}((\diag \boldsymbol{\pi})^{-1} (\C^T)^{-1} \diag(\boldsymbol{\pi}\C) )(\e,\e')=1$
for any $\e\in\EE^*$. 
For diagonal matrices $\mathbf{D}_1$, $\mathbf{D}_2$ and a square matrix $\mathbf{A}$ 
(all of the same sizes) we have   $\mathbf{D}_1 \mathbf{A} \mathbf{D}_2(\e,\e')=$ 
\noindent $\mathbf{D}_1(\e,\e) \mathbf{A}(\e,\e') \mathbf{D}_2(\e',\e')$, thus 
\begin{eqnarray*}
 \sum_{\e'\in\EE^*}((\diag \boldsymbol{\pi})^{-1} (\C^T)^{-1} \diag(\boldsymbol{\pi}\C) )(\e,\e')
 & = & {1\over \pi(\e) }\sum_{\e'\in\EE^*}  \C^{-1}(\e',\e) \diag \boldsymbol({\pi\C})(\e',\e') \\[4pt]
 & = & {1\over \pi(\e) } (\boldsymbol{\pi}\C\C^{-1})(\e)={1\over \pi(\e)}\pi(\e)=1.
\end{eqnarray*}
Thus, $\PP$ is a stochastic matrix and thus $X\sim(\nu,\PP)$ is a Markov chain with the stationary distribution $\pi$. 
Since (\ref{eq:duality}) holds, $X^*$ is an SSD for $X$. 
Theorem \ref{twr:main_twr2_LS} and Corollary \ref{cor:sharp} imply that $X^*$ is a sharp SSD of $X$.
\end{proof}

\begin{rmrk}\label{rem:star_min} \rm
If, in addition, within ordering $\preceq$ we have a unique minimal state, say $\e_1$, and  $X^*$ starts from this state (i.e., $\nu^*=\delta_{\e_1}$),
then the antidual chain also starts from this state, i.e. $\nu=\delta_{\e_1}$. This is the case in all examples that follow.
\end{rmrk}
\begin{rmrk}\label{rem:directlyLink} \rm
The condition that $\widehat{\PP^*}$ is M\"obius monotone (w.r.t. $\pi$ and $\C$) is equivalent to non-negativity
of the resulting matrix $\PP$. In examples, it is often more convenient to calculate $\Lambda$ and $\Lambda^{-1}$ directly.
\end{rmrk}

\subsection{Antidual chains for a generalized coupon collector problem}\label{sec:main_antidual_chains}
Consider $d$ different types of coupons. These are sampled independently with replacement. Sampled types are recorded. For $1\leq k\leq d$ let 
$p_k>0$ be the probability that the coupon of type $k$ is sampled, with $\sum_{k=1}^d p_k\leq 1$. With the remaining probability, i.e., with probability 
$1-\sum_{k=1}^d p_k$, no coupon is sampled.  We start with no coupons of any type.
Let $T^*$ be the number of steps it takes to collect $N_j$ coupons of type $j$, $j=1,\ldots,d$ for some fixed integers $N_1,\ldots,N_d$.
Let $(i_1,\ldots,i_d)$ denote that coupon of type $j$ 
was sampled $i_j$ times.
If $i_j=N_j$ and coupon of type $j$ is sampled, the chain does not move.
The distribution of $T^*$ is the time to absorption in the state  $(N_1,\ldots,N_d)$ of the chain $X^*\sim(\nu^*,\PP^*)$ 
on the state space $\EE^*=\{(i_1,\ldots,i_d): 0\leq i_j\leq N_j, 1\leq j\leq d\}$ with initial distribution
$\nu^*=\delta_{(0,\ldots,0)}$ and the following transition matrix:
\begin{equation}\label{eq:coupon_gen}
\PP^*((i_1,\ldots,i_d),(i'_1,\ldots,i'_d))=
\left\{ 
 \begin{array}{llllllll}
  p_j & \textrm{if} & i_j'=i_j+1, i_k'=i_k, k\neq j, \\[10pt]
  
  \displaystyle 1-\sum_{k=1}^d  p_k + \sum_{k:i_k=N_k} p_k& \textrm{if} & i_j'=i_j, 1\leq j\leq d.
 \end{array}
 \right.
\end{equation}
We refer to $\PP^*$ as to a \textsl{generalized coupon collector chain}.
The case $N_j=1, j=1,\ldots,d$ and $p_k=1/d$ is  \textsl{the classic coupon collector problem}, which has a 
long history, see for example \cite{Feller1971}. 
The term \textsl{generalized} is not unique. It is used when sequence $\{p_k\}$ is general but $N_1=\ldots=N_d=1$ (e.g., \cite{Neal2008}) or when $p_k=1/d$ but we are to collect more 
coupons of each type (see, e.g., \cite{Newman60}, \cite{Doumas2016}).
Although the chain $\PP^*$ given in (\ref{eq:coupon_gen}) includes both mentioned generalizations, we consider 
two antidual chains for two different cases separately: 
\begin{itemize}
 \item[a)]  for  general $N_j\geq 1$ and $p_j, j=1,\ldots,d$  with the uniform stationary distribution of antidual chain;
\item[b)] for general $p_j$ but $N_j=1, j=1,\ldots,d$  with more general stationary distribution of antidual chain (including uniform one as special case).
\end{itemize}

The proofs are postponed to Section \ref{sec:gen_coup_col_proofs}.
\medskip\par 
%
%

For convenience denote $\iii=(i_1,\ldots,i_d)$ and $\iii^{(k)}=(i_1^{(k)},\ldots,i_d^{(k)})$. 
Define $\s_k:=(0,\ldots,1,\ldots,0)$ (with 1 on the position $k$).

\bigskip\par \noindent 
\textbf{Case: general $N_j\geq 1$ and $p_j, j=1,\ldots,d$ and a uniform stationary distribution of   antidual chain}
\medskip\par 

%

\begin{thrm}\label{twr:ver1}
 Let $X^*\sim(\nu^*,\PP^*)$ be a \textsl{generalized coupon collector chain} with the transition matrix given in (\ref{eq:coupon_gen}) with fixed integers $N_j\geq 1, j=1,\ldots,d$.
 Moreover, assume that
\begin{equation}\label{eq:lattice_assumptions}
\sum_{j=1}^d \left(1-{1\over N_j(N_j+1)}\right)p_j\leq 1.
%
\end{equation}
Then the chain $X\sim(\nu,\PP)$ with $\nu=\delta_{(0,\ldots,0)}$ and with transition matrix 
\bigskip\par\noindent
$\PP(\iii^{(1)},\iii^{(2)}) =$
\begin{equation}\label{antidual_gen_NJ}
 \left\{ 
\begin{array}{lllllllllll}
 \displaystyle {i_k^{(1)}+1\over i_k^{(1)}+2}p_k & \textrm{if} & \iii^{(2)}=\iii^{(1)}+\s_k, \\[18pt]
 
 \displaystyle \left({\mathbf{1}(i_k^{(1)}<N_k)\over (i_k^{(1)}+1) (i_k^{(1)}+2)}
 +{\mathbf{1}(i_k^{(1)}=N_k)\over N_k+1}\right)p_k& \textrm{if} & \iii^{(2)}=\iii^{(1)}-m\cdot \s_k   \\
 &  &  \textrm{\ with \ } 1\leq m \leq i_k, \\[12pt]
 
 \displaystyle 1-\sum_{j: i^{(1)}_j<N_j} \left(1-{1\over (i^{(1)}_j+1)(i^{(1)}_j+2)}\right)p_j
- \sum_{j: i^{(1)}_j=N_j}  {N_j\over N_j+1} p_j & \textrm{if} & \iii^{(2)}=\iii^{(1)}
\end{array}\right.
\end{equation}
is an ergodic Markov chain with uniform distribution on $\EE=\EE^*$ which is  a sharp antidual for $\PP^*$.
\end{thrm}
\begin{rmrk}
Note that for example for $N_1=\ldots=N_j=1$, the condition (\ref{eq:lattice_assumptions})
if always  fulfilled.
\end{rmrk}
Roughly speaking, the antidual has the following transitions. Being in state $(i_1,\ldots,i_d)$ it can increase each coordinate by one (if feasible),
it can stay in this state or it can change one of the coordinates to anything smaller. 
Changing some coordinate depends only on the value 
of this coordinate, and decreasing coordinate, say  from $i_j$ to $i_j-m$ is constant for all $1\leq m<i_j$
(the probability depends only on $i_j$ and the formula itself  is different on the border, i.e., when $i_j=N_j$).

\medskip\par \noindent
\textbf{Case: general $p_j$ and $N_j=1, j=1,\ldots,d$ and a non-uniform distribution of antidual chain}.
\medskip\par 

\begin{thrm}\label{twr:coupon_antidual_ver2}
 Let $X^*\sim(\nu^*,\PP^*)$ be a \textsl{generalized coupon collector chain} with the transition matrix given in (\ref{eq:coupon_gen}). Assume that $N_1=\ldots=N_d=1$.
 Let  $a_k\in(0,1)$ for $k=1,\ldots,d$. Then the chain $X\sim(\nu,\PP)$ on the same state space $\EE=\EE^*=\{0,1\}^d$ with the initial distribution
 $\nu=\nu^*=\delta_{(0,\ldots,0)}$ and the transition matrix
\begin{equation}\label{eq:antidual_N1} 
\PP(\iii^{(1)},\iii^{(2)}) =
\left\{ 
\begin{array}{lllllllllll}
 \displaystyle a_k p_k & \rm{if}& \iii^{(2)}=\iii^{(1)}+\s_k, \\[12pt]
 \displaystyle \displaystyle 1-  \sum_{j:i^{(1)}_j=0} a_j p_j -  \sum_{j:i^{(1)}_j = 1} (1-a_j)p_j & \rm{if} & \iii^{(2)}=\iii^{(1)},  \\[12pt]
 \displaystyle (1-a_k) p_k & \rm{if} & \iii^{(2)}=\iii^{(1)}-\s_k, \\
\end{array}\right.
\end{equation}
is an ergodic Markov chain  which is a sharp antidual  for $\PP^*$. The stationary distribution is the following:
\begin{equation}\label{eq:gen_coupon_coll_antidual_pi}
 \begin{array}{lllllll}
\displaystyle \pi(\iii) & = & \displaystyle  \prod_{j=1}^d[a_j \mathbf{1}(i_j=1)+(1-a_j) \mathbf{1}(i_j=0)].\\
\end{array}
\end{equation}
\end{thrm}

\begin{rmrk} \rm
The proof of Theorem \ref{twr:coupon_antidual_ver2} 
implies that the antidual chain $X\sim(\nu,\PP)$ has transitions consistent with partial ordering, i.e.,  at each step it can stay or it can either change one 
coordinate from 0 to 1 or vice-versa. This is not the case for any distribution $\pi$.
It can happen, that for some $\pi$     two coordinates change at a time or antidual 
does not exist (since some entries of $\PP$ can be negative). This is further commented 
after proof in Remark \ref{rm:two_coord}.
\end{rmrk}

\smallskip\par\noindent
Taking the following concrete sequences of $a_k$: $a_k={b\over a+b}$ or $a_k={1\over 2}, j=1,\ldots,d$  we obtain the following special cases:
\begin{crllr}\label{cor:gen_coupon_coll_antidual_special_cases}
The chains $X^{(i)}\sim(\nu,\PP_i), i=1,2$ with a common initial distribution $\nu=\delta_{(0,\ldots,0)}$ and transition matrices

$$
\begin{array}{lllllllll}
 \PP_1(\iii^{(1)},\iii^{(2)})& = &  
 \left\{ 
\begin{array}{llllll}
  \displaystyle{1\over 2} p_k & \mathrm{if} &   \iii^{(2)}=\iii^{(1)}+\s_k, \\[7pt]
  \displaystyle  1-{1\over 2} \sum_{j=1}^d p_j & \mathrm{if} & \iii^{(2)}=\iii^{(1)}, \\[7pt]
  \displaystyle{1\over 2} p_k & \mathrm{if} &  \iii^{(2)}=\iii^{(1)}-\s_k. \\
\end{array}\right.\\
 & \\[10pt]
 \PP_2(\iii^{(1)},\iii^{(2)}) & = & 
\left\{ 
\begin{array}{llllll}
  \displaystyle{b\over {a+b}} p_k & \mathrm{if} &    \iii^{(2)}=\iii^{(1)}+\s_k,\\[7pt]
 \displaystyle 1-{b\over a+b} \sum_{j:i^{(1)}_j=0} p_r - {a\over a+b} \sum_{j:i^{(1)}_j=1} p_r  &  \mathrm{if} & \iii^{(2)}=\iii^{(1)},\\[7pt]
  \displaystyle{a\over {a+b}} p_k & \mathrm{if} &   \iii^{(2)}=\iii^{(1)}-\s_k.\\
\end{array}\right. \\
& \\
\end{array}
$$
and with the  respective stationary distributions
$$ \pi_1(\iii)   =   {1\over 2^d}, \quad \pi_2(\iii) =\displaystyle {a^{d-|\iii|}b^{|\iii|} \over (a+b)^d}$$
%
(where $|\iii|=\sum_{j=1}^d i_j$, called \textsl{a level} of $\iii$) are sharp antidual chains for $\PP^*$ given in (\ref{eq:coupon_gen}).
\end{crllr}

\begin{rmrk}\rm
In   \cite{Lorek2012d} we considered the chain on $\EE=\{0,1\}^d$ with   transition matrix $\PP_3$ given by
 $$
 \PP_3(\iii^{(1)},\iii^{(2)})=\left\{ 
\begin{array}{llllll}
  \displaystyle \alpha_k & \mathrm{if} &   \iii^{(2)}=\iii^{(1)}+\s_k, \\[7pt]
  \displaystyle 1-\sum_{j: i^{(1)}_j=0} \alpha_j - \sum_{j: i^{(1)}_j=1} \beta_j    & \mathrm{if} & \iii^{(2)}=\iii^{(1)}, \\[7pt]
  \displaystyle \beta_k & \mathrm{if} &  \iii^{(2)}=\iii^{(1)}-\s_k. \\
\end{array}\right.
$$
The chain is reversible with product form stationary distribution:
\begin{equation*}\label{eq:nonsym_cube_stat}
\pi_3(\iii)=  \prod_{j: i_j=0}{\beta_j\over \alpha_j+\beta_j}\prod_{j: i_j=1}{\alpha_j\over \alpha_j+\beta_j}.
\end{equation*}
We showed that the chain is M\"obius monotone if and only if $\sum_{j=1}^d(\alpha_j+\beta_j) \leq 1$. As partial ordering,
coordinate-wise was used. Then we obtained the following dual chain:
 $$
 \PP^*(\iii^{(1)},\iii^{(2)})=\left\{ 
\begin{array}{llllll}
  \displaystyle \alpha_k+\beta_k & \mathrm{if} &   \iii^{(2)}=\iii^{(1)}+\s_k, \\[7pt]
  \displaystyle 1-\sum_{j: i^{(1)}_j=0} (\alpha_j+\beta_j)   & \mathrm{if} & \iii^{(2)}=\iii^{(1)},  
\end{array}\right.
$$
what is our absorbing dual (\ref{eq:coupon_gen}) we started with, with $p_j=\alpha_j+\beta_j$ and $N_j=1, j=1,\ldots,d$.
Note that $\PP_3$  is a special case of $\PP$ given in (\ref{eq:antidual_N1}) with $a_j={\alpha_j\over \alpha_j+\beta_j}$.
%
\end{rmrk}
%
%
%
%
%
\begin{crllr}\label{cor:P123_eigenv}
 The matrices $\PP$ given in (\ref{antidual_gen_NJ})  and in (\ref{eq:antidual_N1}) have  eigenvalues  of the form:
 $$ \lambda_A=1-\sum_{k\in A} p_k, \quad \textrm{ for } \ A\subseteq\{1,\ldots,d\}$$
 (the multiplicity of which depends on the case).
\end{crllr}
\begin{proof}
 We can order the states of $X^*$ in such a way  that $\PP^*$ given in (\ref{eq:coupon_gen}) is upper triangular, thus
 the
eigenvalues are the entries on the diagonal. 
If the link $\Lambda$ is invertible (which is the case), then the transition matrices $\PP$ and  $\PP^*$ of SSD have the same set of eigenvalues,
what is a direct consequence of relation (\ref{eq:duality}).
\end{proof}

\begin{rmrk}
 \rm Fix $d$ and $N_j=N, j=1,\ldots,d$. One can ask the following question: 
 For what sequence $\{p_k\}$ is the associated $T_{FSST}$  stochastically the smallest? Conjecture 2 in \cite{Doumas2016} suggests that 
 this is in the case of equal probabilities $p_k=1/d$.
\end{rmrk}

\subsection{Results on the separation cutoff}\label{sec:cutoff_results}
Since obtained antidual chains are sharp (i.e., (\ref{eq:fsst}) holds), we can present a series of results on the separation cutoff utilizing existing results 
on the limiting distribution of $T^*$. 
\par 
We start with the simplest chain corresponding to the classical coupon collector problem.
\begin{crllr}\label{cor:classic_antidual}
Consider a sequence of Markov chains $X_{(d)}$ indexed by $d=1,2,\ldots$ on 
$\EE_{(d)}=\{0,1\}^d$ with an initial distribution $\nu_{(d)}=\delta_{(0,\ldots,0)}$ 
and the transition matrix $\PP_{(d)}$ given in (\ref{eq:antidual_N1}) with $p_k={1\over d}$ and any $a_k\in(0,1)$ for $k=1,\ldots, d$.
The stationary distribution $\pi_{(d)}$ is given in (\ref{eq:gen_coupon_coll_antidual_pi}).
The sequence exhibits a separation cutoff at time $d\log d$ with window size $d$.
\end{crllr}
\begin{proof}
Denote the FSST of the chain by $T^*_d$. It is known that $ET^*_d=d\sum_{i=1}^d{1\over i}\approx d\log d$. Moreover, 
${1\over d}(T^*_d-d\log d)$ converges in distribution (as $d\to\infty$) to  a standard Gumbel random variable $Z$ (with c.d.f $P(Z\leq c)=e^{-e^{-c}}$),
see \cite{Holst2001}.\par 
\noindent
Taking $t_d=d\log d$ and $w_d=d$ we have
$$
\begin{array}{llllllllll}
 sep(\nu_{(d)}\PP_d^{d\log d+c d},\pi_d) & = & P(T_d^*>d\log d + cd)=1-P\left({1\over d}(T_d^*-d\log d) \leq  c\right), \\[10pt]
sep(\nu_{(d)}\PP_d^{d\log d-c d},\pi_d) & = & P(T_d^*>d\log d - cd)= 1-P\left({1\over d}(T_d^*-d\log d) ) \leq -c\right).
\end{array}
$$
Taking the limits as $d\to\infty$ we have
$$
\begin{array}{llllll}
 \displaystyle   \limsup_{d\to\infty} & sep(\nu_{(d)}\PP_d^{d\log d + c d},\pi_d)  & = & 1-e^{-e^{-c}},\\[8pt]
 \displaystyle \liminf_{d\to\infty} & sep(\nu_{(d)}\PP_d^{d\log d-  c d},\pi_d)  & = & 1-e^{-e^{c}}.
 
\end{array}
$$
Taking the limit as $c\to\infty$ finishes the proof.
\end{proof}

Results on the limiting distribution of $T^*_d$ from \cite{Neal2008} let us indicate separation cutoffs for cases with non-constant probabilities $p_k$.
For example we can have the following corollary.
\begin{crllr}
Consider  piecewise constant probability density function on $[0,1]$:
$$f(y)=\lambda_j, \quad n_{j-1}<x\leq n_j, \quad 1\leq j\leq k,$$
where $\lambda_1,\ldots,\lambda_k>0$ and $0=n_0<n_1<\cdots n_k=1$. Without loss of generality assume that $\lambda_1 <\lambda_2<\ldots<\lambda_k$. 
Consider a sequence of Markov chains $X_{(d)}$ indexed by $d=1,2,\ldots$ on
$\EE_{(d)}=\{0,1\}^d$ with an initial distribution $\nu_{(d)}=\delta_{(0,\ldots,0)}$ 
and the transition matrix $\PP_{(d)}$ given in (\ref{eq:antidual_N1}) with 
$$p_k=\int_{(k-1)/d}^{k/d} f(y)dy, \quad k=1,\ldots, d$$
and any $a_k\in(0,1)$ for $k=1,\ldots, d$. The stationary distribution $\pi_{(d)}$ is given in (\ref{eq:gen_coupon_coll_antidual_pi}).
The sequence exhibits a separation cutoff at time $t_d={d\over \lambda_1}(\log d-\log(n_1))$ with window size  $w_d={d\over\lambda_1}$.
\end{crllr}
\begin{proof}
Denote the FSST of the chain by $T^*_d$ (which is equal, in distribution, to collecting $d$ coupons).
We have 
$$
\begin{array}{llllllllll}
  sep(\nu_{(d)}\PP_d^{t_d+cw_d},\pi_d)  & = & P\left(T_d^*>{d\over \lambda_1}(\log d-\log(n_1))+c{d\over\lambda_1}\right) \\[10pt]
 & = & 1-P\left({1\over d}(T_d^*-{1\over \lambda_1}d\log d) \leq {\log(n_1)\over \lambda_1} + {c\over\lambda_1}\right).
\end{array}
$$
Lemma 3.1 in \cite{Neal2008} 
implies that ${1\over d} (T^*_d-{1\over \lambda_1}d\log d)$  converges in 
distribution to a random variable $Z$ with c.d.f
$P(Z\leq c) =e^{-{n_1 e^{-\lambda_1 c}}}$.
Thus, we have 
$$\limsup_{d\to\infty}   sep(\nu_{(d)}\PP_d^{t_d+cw_d},\pi_d) 
=1-e^{-{n_1 e^{-\lambda_1 \left({\log(n_1)\over \lambda_1} + {c\over\lambda_1}\right)}}}=1-e^{-e^{-c}}.$$
Similarly
$$ \textstyle sep(\nu_{(d)}\PP_d^{t_d-cw_d},\pi_d)= 1-P\left({1\over d}(T_d^*-{1\over \lambda_1}d\log d) \leq {\log(n_1)\over \lambda_1} - {c\over\lambda_1}\right)
$$
and
$$\liminf_{d\to\infty}   sep(\nu_{(d)}\PP_d^{t_d-cw_d},\pi_d) 
=1-e^{-{n_1 e^{-\lambda_1 \left({\log(n_1)\over \lambda_1} - {c\over\lambda_1}\right)}}}=1-e^{-e^{c}}.$$
Taking limits as $c\to\infty$ finishes the proof.
\end{proof}
\par\noindent 
Next corollaries utilize    results on time until some set of coupons is collected.

\begin{crllr}\label{cor:antidual_Ncoupons_symm}
Consider a sequence of Markov chains $X_{(d)}$ indexed by $d=1,2,\ldots$ on $\EE_{(d)}=\{0,1,\ldots,N\}^d$ with 
an initial distribution $\nu_{(d)}=\delta_{(0,\ldots,0)}$ 
and the transition matrix $\PP_{(d)}$ given in (\ref{antidual_gen_NJ}) with 
$p_k={1\over d}$ and $N_1=\ldots=N_d=N\geq 2$ (so that (\ref{eq:lattice_assumptions}) holds).
The stationary distribution $\pi_{(d)}$ is uniform. 
The sequence of chains exhibits a separation cutoff at time 
$d\log d +(N-1)d\log\log d $ with window size $d$.
\end{crllr}
\begin{proof}
In \cite{Erdos1961} authors derived limiting distribution of $T^*_d$ showing that 
$$ {1\over d} ({T^*_d- d\log d - (N-1)d\log\log d + d[\gamma-\log(N-1)!]})$$ 
(where $\gamma=0.57721\ldots$ is the Euler-Mascheroni constant)
converges in distribution to a standard Gumbel random variable.
Similar calculations as in Corollary \ref{cor:classic_antidual} finish the proof.

\end{proof}

Recently authors in \cite{Doumas2016} extended the result of \cite{Erdos1961}
obtaining the limiting distribution of $T^*_d$ for $N_1=\ldots=N_d=N$ and for 
quite general choices of probabilities $p_k$. Let us  indicate here one example
(which actually includes result of Corollary \ref{cor:antidual_Ncoupons_symm} as a special case).

\begin{crllr}\label{cor:antidual_Ncoupons_nons_log}
Consider a sequence of Markov chains $X_{(d)}$ indexed by $d=1,2,\ldots$ on $\EE_{(d)}=\{0,1,\ldots,N\}^d$ with 
an initial distribution $\nu_{(d)}=\delta_{(0,\ldots,0)}$ 
and the transition matrix $\PP_{(d)}$ given in (\ref{antidual_gen_NJ}) with 
$$p_k={1\over (\log k)^p} {1\over K_d}, \quad K_d=\sum_{k=1}^d {1\over (\log k)^p}, \quad p\in(0,1), \quad k=1,\ldots,d$$
 and $N_1=\ldots=N_d=N\geq 2$ (so that (\ref{eq:lattice_assumptions}) holds).
The stationary distribution $\pi_{(d)}$ is uniform. 
The sequence of chains exhibits a separation cutoff at time 
$d\log d +(N-1)d\log\log d $ 
with window size $d$.
\end{crllr}
\begin{proof}
In \cite{Doumas2016} authors prove that 
$$ {1\over d} ({T^*_d- d\log d - (N-1)d\log\log d + d[\gamma+p-\ln(p+1)-\log(N-1)!]})$$ 
converges in distribution to a standard Gumbel random variable.
Again, similar calculations as in Corollary \ref{cor:classic_antidual} finish the proof.
\end{proof}

\subsection{Constructing an ergodic chain with a prespecified FSST and an arbitrary stationary distribution}\label{sec:constr_FSST}
Let us ask the following question 
(which was one of the main motivations for the paper):
 \begin{center}
  \begin{minipage}[t]{0.9\textwidth}
\textsl{How to construct a Markov chain 
on a state space of size $M$ with arbitrary stationary distribution $\pi$ whose FSST $T$ is deterministic, $P(T=M-1)=1$?}
\end{minipage}  
 \end{center}
 
 The recipe is clear from previous sections: Start with some absorbing chain $X^*$ for which $P(T^*=M-1)=1$, where $T^*$ 
 is the absorption time. Probably the simplest one is  the following: take $\EE=\{1,\ldots,M\}$ with transitions 
 $\PP_0^*(k,k+1)=1$ for $k<N$ and $\PP_0^*(N,N)=1$ and start it at state 1.  Then of course we have desired 
 absorption time and thus the antidual would have desired stationary distribution and FSST.
 \par 
 The above example will be a special case of a more general result. Many absorbing chains have
 the  absorption time $T^*$ distributed as a mixture of sums of independent geometric 
 random variables with parameters being the eigenvalues of the transition matrix.
 E.g., for stochastically monotone discrete time birth and death chain starting at 1 with $d>1$ being the 
 absorbing state, the
 time to absorption is distributed as a sum of geometric random variables with parameters 
 being the eigenvalues of the transition matrix (which are positive in this case). This result follows from  Karlin and McGregor \cite{Karlin1959} or
 Keilson \cite{Keilson71b}.
 Fill \cite{Fill2009} gave a first stochastic proof of this result using dualities
 (the result was simultaneously obtained in \cite{Diaconis2009b}).
 This was extended to skip-free Markov chains in Fill \cite{Fill2009a}.
 Miclo \cite{Miclo2010a} showed that for any absorbing chain on $\EE=\{\e_1,\ldots,\e_M\}$
 with  positive eigenvalues and some reversibility condition (involving substochastic kernel
 corresponding to the transition matrix with row and column corresponding to absorbing state removed) there exists a measure 
 $a=(a_1,\ldots,a_M)$ such that the time to absorption $T^*$ has distribution 
 $$T^*\sim \sum_{i=1}^{M-1} a_i \mathcal{G}(\lambda_i,\lambda_{i+1},\ldots,\lambda_{M-1}),$$
 where $\lambda_i$ are the eigenvalues of the transition matrix  sorted in non-increasing order and $\mathcal{G}(p_1,\ldots,p_k)$ denotes 
 the distribution of $\sum_{j=1}^k X_j, $ where $X_j \sim Geo(p_j)$.
 \bigskip\par \noindent
 For convenience denote $H(k):=\sum_{j=1}^k \pi(j)$. Our result is following.
 \begin{thrm}\label{thm:FSST_gen}
  Let $\EE=\{1,\ldots,M\}$ and $p_k\in[0,1], k=1,\ldots,{M-1}$. Let $a_k, \pi(k), k=1,\ldots,M$ be two probability
  distributions on $\EE$ such that $a_k\geq 0, \pi(k)>0$ for all $k\in\EE$. 
  Define the matrix 
    $$ \PP(k,s) = 
  \left\{ 
  \begin{array}{llllllll}
  \displaystyle 1-{\pi(2)\over \pi(1)+\pi(2)}p_1 & \textrm{if} & k=s=1, \\[12pt]
  \displaystyle {\pi(2)\over \pi(1)+\pi(2)}p_1 & \textrm{if} & k=1, s=2, \\[12pt]
  \displaystyle {\pi(s)\over \pi(k)}\left[ p_{k-1}\left( 1-{H(k-1)\over H(k)}\right)-p_{k} \left(1-{H(k)\over H(k+1)}\right) \right] & \textrm{if} & 1<k<M, s<k, \\[12pt]
  \displaystyle 1-p_k\left(1 - {H(k)\over H(k+1)}\right)-p_{k-1} {H(k-1)\over H(k)} & \textrm{if} & 1<k<M, s=k, \\[12pt]
   \displaystyle p_k {H(k)\over  H(k+1)}   {\pi(k+1)\over \pi(k)}   & \textrm{if} & 1<k<M, s=k+1, \\[12pt]
 \displaystyle     p_{M-1}\pi(s) & \textrm{ if } & k=M, s\leq M-1, \\[12pt]
 \displaystyle    1-p_{M-1}+p_{M-1}\pi(M)  & \textrm{ if } & k=M, s=M.         
  \end{array}
\right. 
$$
Assume that $\pi$ and sequence $\{p_k\}_{k=1,\ldots,M}$ are such that that the matrix $\PP$ is non-negative. Then Markov 
chain $X$ with the  transition matrix $\PP$ and with  the initial distribution $\nu=(\nu(1),\ldots,\nu(M))$ given by
$$\nu(k)=\pi(k)\sum_{i=k}^M  {a_i\over H(i)}$$
has the FSST $T$ distributed as 
\begin{equation}\label{eq:FSST_T}
\sum_{i=1}^{M-1} a_i \mathcal{G}(p_i,p_{i+1},\ldots,p_{M-1})
\end{equation}
and $\pi$ is its stationary distribution. Moreover, $\{1-p_1,\ldots,1-p_{M-1}, 1\}$ are the eigenvalues of $\PP$.
 \end{thrm}
  Note that $X$ is a skip-free chain: for given $k$ the only nonzero entries of $\PP$ are $\PP(k,s)$ for $s\leq k+1$.
  The proof of the theorem is postponed to Section \ref {sec:proof_FSST}. \par 
  \noindent We can relatively easy have  some corollaries being interesting special cases of Theorem \ref{thm:FSST_gen}.
  Applying the Theorem \ref{thm:FSST_gen} with $p_k=1, k=1,\ldots,M-1, p_M=0$ and $a_1=1, a_k=0, k=2,\ldots,M$ we obtain the following corollary.
  \begin{crllr}\label{cor:pure}
 Consider a distribution $\pi$  on $\EE=\{1,\ldots,M\}$ such that $\pi(k)>0$ for all $k\in\EE$.
 The  Markov chain $X$ on $\EE$ with transition matrix 
 $$\PP_0(k,r)=
 \left\{
 \begin{array}{llll}
 \displaystyle {\pi(r)\over \pi(1)+\pi(2)} & \textrm{ for } & k=1, r\in\{1,2\},\\[12pt]
  \displaystyle{\pi(r)\over \pi(k)}\left[{H(k)\over H(k+1)}-{H(k-1)\over H(k)}\right] & \textrm{ for } & 1<k<M, r\leq k, \\[12pt]
  \displaystyle{\pi(k+1)\over \pi(k)}{H(k)\over H(k+1)} & \textrm{ for } & 1<k<M, r=k+1, \\[12pt]
  \displaystyle\pi(r) & \textrm{ for } & r\leq k=M \\
 \end{array}\right.
$$ 
is ergodic with  the stationary distribution $\pi$. 
Assume the initial distribution is $\nu=\delta_1$ (i.e., $P(X_0=1)=1$). Then the chain has deterministic fastest strong stationary time
$T$ such that $P(T=M-1)=1$.
\end{crllr}

Note that for this chain we have 
$$sep(\nu\PP^k,\pi)=P(T>k)
=\left\{\begin{array}{llll}
  1 & \textrm{ if } k\leq M-2, \\
  0 & \textrm{ if } k\geq M-1. \\
 \end{array}\right.
$$
Thus, this is an extreme example for a separation cutoff: For any $k\leq M-2$ the chain is completely not mixed 
(the separation   between stationary distribution and distribution at step $k$ is 1) and the chain mixes completely exactly 
at step $k=M-1$ (the distance is 0).
\par 
\noindent 
Simplifying the chain further by taking additionally uniform distribution $\pi(k)={1\over M}$ in Corollary \ref{cor:pure} we obtain

 $$\PP_0(k,r)=
 \left\{
 \begin{array}{llll}
 \displaystyle {1\over 2} & \textrm{ for } & k=1, r\in\{1,2\},\\[12pt]
  \displaystyle{1\over k(k+1)} & \textrm{ for } & 1<k<M, r\leq k, \\[12pt]
  \displaystyle{k\over k+1} & \textrm{ for } & 1<k<M, r=k+1, \\[12pt]
  \displaystyle{1\over M} & \textrm{ for } & r\leq k=M. \\
 \end{array}\right.
$$
The chain is sketched in Fig. \ref{fig:uniform_x}
\begin{center}
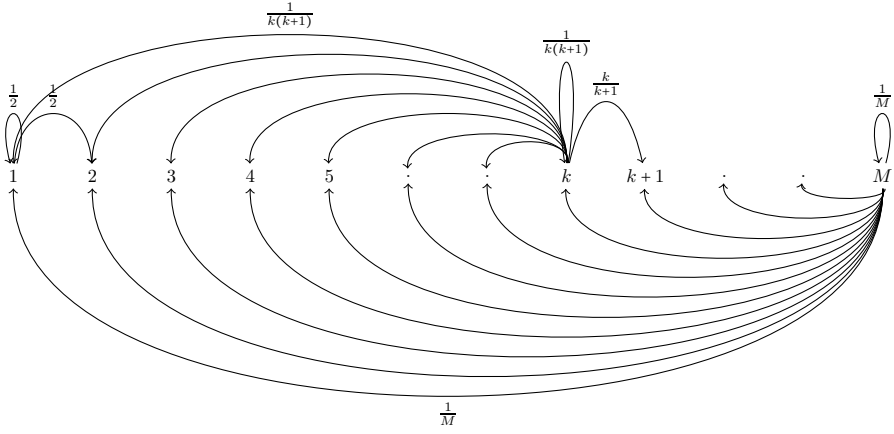
\begin{figure}\label{fig:gambler1d}
\begin{tikzpicture}[line join=round,x=1.0cm,y=1.0cm,scale=0.7,every node/.style={scale=0.7}]
 \node (p1)  at (-2,0cm) {$1$};
 \node (p2) at (-0.5,0cm) {$2$};
 \node (p3) at (1,0cm) {$3$};
 \node (p4) at (2.5,0cm) {$4$};
 \node (p5) at (4.0,0cm) {$5$};
 \node (p6) at (5.5,0cm) {$.$};
 \node (p7) at (7.0,0cm) {$.$};
 \node (pk) at (8.5,0cm) {$k$};
 \node (pk1) at (10.0,0cm) {$k+1$};
 \node (pk2) at (11.5,0cm) {$.$};
 \node (pk3) at (13.0,0cm) {$.$};
 \node (pM) at (14.5,0cm) {$M$};

  \draw [->] (p1) .. controls +(0.4,1.5) and +(-0.4,1.5) .. node [midway, above] {${1\over 2}$} (p1);
  \draw [->] (p1) .. controls +(0.1,1.5) and +(-0.1,1.5) .. node [midway, above] {${1\over 2}$} (p2);
  \draw [->] (pk) .. controls +(0.1,3.5) and +(-0.1,3.5) .. node [midway, above] {${1\over k(k+1)}$} (p1);
  \draw [->] (pk) .. controls +(0.1,3.0) and +(-0.1,3.0) .. node [midway, above] {} (p2);
  \draw [->] (pk) .. controls +(0.1,2.5) and +(-0.1,2.5) .. node [midway, above] {} (p3);
  \draw [->] (pk) .. controls +(0.1,2.0) and +(-0.1,2.0) .. node [midway, above] {} (p4);
  \draw [->] (pk) .. controls +(0.1,1.5) and +(-0.1,1.5) .. node [midway, above] {} (p5);
  \draw [->] (pk) .. controls +(0.1,1.0) and +(-0.1,1.0) .. node [midway, above] {} (p6);
  \draw [->] (pk) .. controls +(0.1,0.8) and +(-0.1,0.8) .. node [midway, above] {} (p7);
  \draw [->] (pk) .. controls +(0.4,2.8) and +(-0.4,2.8) .. node [midway, above] {${1\over k(k+1)}$} (pk);
  \draw [->] (pk) .. controls +(0.4,1.8) and +(-0.4,1.8) .. node [midway, above] {${k\over k+1}$} (pk1);
  \draw [->] (pM) .. controls +(0.1,-5.5) and +(-0.1,-5.5) .. node [midway, below] {${1\over M}$} (p1);
  \draw [->] (pM) .. controls +(0.1,-5) and +(-0.1,-5) .. node [midway, above] {} (p2);
  \draw [->] (pM) .. controls +(0.1,-4.5) and +(-0.1,-4.5) .. node [midway, above] {} (p3);
  \draw [->] (pM) .. controls +(0.1,-4) and +(-0.1,-4) .. node [midway, above] {} (p4);
  \draw [->] (pM) .. controls +(0.1,-3.5) and +(-0.1,-3.5) .. node [midway, above] {} (p5);
  \draw [->] (pM) .. controls +(0.1,-3.0) and +(-0.1,-3.0) .. node [midway, above] {} (p6);
  \draw [->] (pM) .. controls +(0.1,-2.5) and +(-0.1,-2.5) .. node [midway, above] {} (p7);
  \draw [->] (pM) .. controls +(0.1,-2.0) and +(-0.1,-2.0) .. node [midway, above] {} (pk);
  \draw [->] (pM) .. controls +(0.1,-1.5) and +(-0.1,-1.5) .. node [midway, above] {} (pk1);
  \draw [->] (pM) .. controls +(0.1,-1.0) and +(-0.1,-1.0) .. node [midway, above] {} (pk2);
  \draw [->] (pM) .. controls +(0.1,-0.5) and +(-0.1,-0.5) .. node [midway, above] {} (pk3);
  
  \draw [->] (pM) .. controls +(0.4,1.5) and +(-0.4,1.5) .. node [midway, above] {${1\over M}$} (pM);
  
\end{tikzpicture}
\caption{The chain $X$ on $\EE=\{1,\ldots,M\}$ with the uniform stationary 
distribution and with the deterministic FSST $T: P(T=M-1)=1$}\label{fig:uniform_x}
\end{figure}
\end{center}
\smallskip\par 
\noindent\textbf{Two Markov chains on essentially different state spaces with the same FSST} \
\par 
\noindent 
So far in this section we considered chains on totally ordered state space $\EE=\{1,\ldots,M\}$. We can also consider another state spaces.
We will consider chain on $\EE^{(2)}=\{0,1\}^d$. We will not present full generality one can have,
instead we will present two chains, one on $\EE^{(1)}=\{1,\ldots,d\}$ and the other on $\EE^{(2)}$
both with  uniform distributions and the same FSST distributed as $\sum_{k=1}^{d-1} X_k, $ where $X_k\sim Geo(k\cdot p)$ 
for some fixed $p\leq {1\over d}$. Note that in particular the sizes of the state spaces are completely
different, $2^d$ versus $d$\par 
\begin{crllr}
Fix some integer $d>1$ and $0<p\leq {1\over d}$.
 Let  $X^{(1)}$ be a Markov chain on $\EE^{(1)}=\{1,\ldots,d\}$ with an initial distribution $\nu^{(1)}=(1,0,\ldots,0)$ and transitions
       $$ \PP^{(1)}(k,s) = 
  \left\{ 
  \begin{array}{llllllll}
  \displaystyle 1-{1\over 2}(d-1)p & \textrm{if} & k=s=1, \\[12pt]
  \displaystyle {1\over 2}(d-1)p & \textrm{if} & k=1, s=2, \\[12pt]
  \displaystyle { k(d-k+1)+1\over k(k+1)}p& \textrm{if} & 1<k<d, s<k, \\[12pt]
  \displaystyle 1-\left({d-k\over k+1} + {(d-k+1)(k-1)\over k}\right) p   & \textrm{if} & 1<k<d, s=k, \\[12pt]
   \displaystyle (d-k)p{k\over k+1}    & \textrm{if} & 1<k<d, s=k+1, \\[12pt]
 \displaystyle   {p\over d} & \textrm{ if } & k=d, s\leq d-1, \\[12pt]
 \displaystyle   {p\over d}+1-p   & \textrm{ if } & k=d, s=d.         
  \end{array}
\right. 
$$
Let  $X^{(2)}$ be a Markov chain on $\EE^{(2)}=\{0,1\}^d$  with initial distribution 
 $\nu^{(2)}((0,\ldots,0))=\nu^{(2)}((1,0,\ldots,0))=1/2$ and with transitions
$$ \PP^{(2)}(\iii,\iii') = 
  \left\{\begin{array}{llllllll}
  {1\over 2}p & \textrm{if} & \iii'=\iii \pm \mathbf{s}_k, \\[10pt]
   1-{1\over 2} dp  &\textrm{if} & \iii'=\iii.
\end{array}\right. 
$$
(Recall that $|\iii|=\sum_{j=1}^d i_j$   was called a level of $\iii$). 
\smallskip\par \noindent
Then the FSSTs $T^{(1)}$ and $T^{(2)}$ of both chains have the same distribution:
 $$T^{(1)} \stackrel{(d)}{=}T^{(2)} = \sum_{k=1}^{d-1} X_k, \ \textrm{where } \ X_k\sim Geo(k\cdot p).$$
 Both chains have the uniform stationary distribution on respective state spaces.
\end{crllr}
\begin{proof}
\noindent
We will show that  chains $X^{(1)}$ and $X^{(2)}$ are sharp antidual chains of different chains $X^{*(1)}$ and $X^{*(2)}$,
whose absorption times are equal to the statement. 
\begin{itemize}
 \item Chain $X^{(1)}$
 
This is a special case  of the chain given in Theorem \ref{thm:FSST_gen}
with $p_k=(d-k)p$ and the uniform stationary distribution $\pi$.
Taking $a_1=1, a_k=0, k=2,\ldots,M$ we have that the initial distribution $v=(1,0,\ldots,0)$ and that 
FSST $T^{(1)}$ is distributed as $\sum_{k=1}^{d-1} X_k, X_k\sim Geo(p_k)$ with $p_k=(d-k)p$. 
The distribution of $T^{(1)}$ is equal to $\sum_{k=1}^{d-1} Y_k$ with $Y_k\sim Geo(k\cdot p)$
 \item Chain $X^{(2)}$

This is a special case of the chain $\PP_1$ given in 
Corollary \ref{cor:gen_coupon_coll_antidual_special_cases} with $p_k=p$. Thus, its sharp dual chain is given in 
(\ref{eq:coupon_gen}). Recall this is the case $N_j=1, j=1,\ldots,d$, let us explicitly write the transitions 
of this $\PP^*$ using notation from this section:
$$\PP^*(\iii,\iii')=
\left\{ 
 \begin{array}{llllllll}
  p & \textrm{if} & \iii'=\iii + \mathbf{s}_k, \\[10pt]
  1-(d-|\iii|)p  &\textrm{if} & \iii'=\iii
 \end{array}
 \right.
 $$
 Roughly speaking, this is the following random walk on hypercube $\{0,1\}^d$. Being at some state 
 $\iii=(i_1,\ldots,i_d), i_k\in\{0,1\}$ either we change one coordinate from 0 to 1 with probability 
 $p$ or with the remaining probability we do nothing. State $(1,\ldots,1)$ is an absorbing state. 
 Since the probability of changing 0 into 1 does not depend on the actual state, the
 time to increase the current level depends only on the level. Being at any state on level $|\iii|=l$
 the time to reach next level has distribution $Geo((d-l)p)$ (since there are $(d-l)$ of zeros, each 
 of which can be changed into 1 with probability $p$). Thus, if the chain starts 
 somewhere on level 1, say $\nu^*((1,0,\ldots,0))=1$, then  the absorption time  is equal 
 in distribution to $\sum_{k=1}^{d-1} X_k,$ where $X_k \sim Geo(k\cdot p)$. 
 What remains to show is that 
 $\nu=\nu^*\Lambda$ yields $\nu^{(2)}((0,\ldots,0))=\nu^{(2)}((1,0,\ldots,0))=1/2$. 
 All the proofs of Theorems    \ref{twr:ver1} and \ref{twr:coupon_antidual_ver2} are based on coordinate-wise ordering, i.e.,
  \begin{equation}\label{eq:coord_ordering}
 \iii\preceq\iii' \textrm{ if  } i_j\leq i'_j, j=1,\ldots,d. 
 \end{equation}
 Recall the link $\Lambda$ (it is given in (\ref{eq:link0}))
$$\Lambda(\iii,\iii')= {\pi(\iii')\over \sum_{\iii_0:\iii_0\preceq \iii} \pi(\iii_0)}\mathbf{1}(\iii'\preceq \iii).$$
We have 
 $$
 \begin{array}{lllll}
  \nu(0,\ldots,0)&=&\displaystyle \sum_{\iii} \nu^*(\iii)\Lambda(\iii,(0,\ldots,0)) = \Lambda((1,0,\ldots,0),(0,\ldots,0)) \\[10pt]
		&=& \displaystyle  {\pi((0,\ldots,0))\over \pi((0,\ldots,0))+\pi((1,0,\ldots,0))}={1\over 2},\\[20pt]
 \nu(1,0,\ldots,0)& = & \displaystyle\sum_{\iii} \nu^*(\iii)\Lambda(\iii,(1,0,\ldots,0))=\Lambda((1,0,\ldots,0),(1,0,\ldots,0))\\[10pt] 
		    &=&\displaystyle {\pi((1,0,\ldots,0))\over \pi((0,\ldots,0))+\pi((1,0,\ldots,0))}={1\over 2},\\ 		
 \end{array}
$$
what finishes the proof.

\end{itemize}

\end{proof}

\section{Proofs}\label{sec:proofs}
\subsection{Proofs of Theorems  \ref{twr:ver1} and \ref{twr:coupon_antidual_ver2}}\label{sec:gen_coup_col_proofs}
\noindent 
In both proofs we use the coordinate-wise ordering (defined in (\ref{eq:coord_ordering})) for which   $\iii_{min}=(0,\ldots,0)$ is the
unique minimal and $\iii_{max}=(N_1,\ldots,N_d)$ is the unique maximal one.
\begin{proof}[Proof of Theorem \ref{twr:ver1}]
For the ordering under consideration, directly from Proposition 5 in Rota \cite{Rota64}, we find the
corresponding M\"obius function
\begin{equation}\label{eq:coord_wise_mobius_fun}
\mu((i_1,\ldots,i_d), (i_1+r_1,\ldots,i_d+r_d))
=\left\{\begin{array}{ll}
  \displaystyle(-1)^{\sum_{k=1}^d r_k}  &  r_j\in\{0,1\}, \ i_j+r_j\leq N_j,\  k=1,\ldots,d \\[8pt]
0 & \textrm{otherwise}. \\
\end{array}\right. 
\end{equation}
Let
$$\rho(\iii)=\prod_{j=1}^d (i_j+1).$$
 We will apply Theorem \ref{twr:main} with the above ordering and the uniform distribution $\pi$ on 
$\EE^*$, i.e., $\pi(\iii)={1\over \rho(\iii_{max})}$. Since $X^*$ starts at the minimal state, so does - by Remark \ref{rem:star_min} - the antidual chain.
\noindent
The link $\Lambda(\iii,\iii')$
is the uniform distribution truncated to $\{\iii'\preceq \iii \}$, from (\ref{eq:link}) we have
$\Lambda=(\diag(\pi\C))^{-1} \C^T \diag(\pi)$, thus 
$$\Lambda(\iii,\iii')={\C(\iii',\iii) {1\over \rho(\iii_{max})}\over \displaystyle \sum_{\iii^{(2)}} {1\over \rho(\iii_{max})} \C(\iii^{(2)},\iii)}={\mathbf{1}(\iii'\preceq \iii)\over \rho(\iii)}.$$
The inverse is given by $\Lambda^{-1}=(\diag(\pi))^{-1}  (\C^{-1})^T \diag(\pi\C)$, thus 
$$\Lambda^{-1}(\iii^{(1)},\iii^{(2)})={1\over {1\over  \rho(\iii_{max})}}\C^{-1}(\iii^{(2)},\iii^{(1)}) {\rho(\iii^{(2)})\over \rho(\iii_{max})}= \rho(\iii^{(2)}) \C^{-1}(\iii^{(2)},\iii^{(1)}).$$
Instead of calculating $\widehat{\PP^*}$, we will calculate $\Lambda^{-1}$ and then 
directly the antidual chain  from $\PP=\Lambda^{-1}\PP^*\Lambda$
(the conditions on $(\pi,\C)$-M\"obius monotonicity will be read from the resulting antidual, see Remark \ref{rem:directlyLink}).
We have to calculate
$$\PP(\iii^{(1)},\iii^{(2)})=(\Lambda^{-1}\PP^*\Lambda  )(\iii^{(1)},\iii^{(2)})=\sum_{\iii} \Lambda^{-1}(\iii^{(1)},\iii) (\PP^*\Lambda ) (\iii,\iii^{(2)}).$$
Because of the form of $\Lambda^{-1}$, we need only to consider states which differ from $\iii^{(1)}$ at most by 1 on each coordinate.

 \smallskip\par 
 
 $$
 \begin{array}{lrll} 
\PP(\iii^{(1)},\iii^{(2)})   & = & \displaystyle\sum_{\rr=(r_1,\ldots,r_d)\in\{0,1\}^d:\ \iii^{(1)}-\rr\in\EE^*  } \Lambda^{-1}(\iii^{(1)},\iii^{(1)}-\rr) (\PP^*\Lambda ) (\iii^{(1)}-\rr,\iii^{(2)}) \\[18pt]
 & = & \displaystyle\sum_{\rr=(r_1,\ldots,r_d)\in\{0,1\}^d:\ \iii^{(1)}-\rr\in\EE^* } (-1)^{|\rr|}\rho(\iii^{(1)}-\rr) (\PP^*\Lambda ) (\iii^{(1)}-\rr,\iii^{(2)}). \\[8pt]
\end{array}
$$

\noindent 
 We need to calculate 
 $$ (\PP^*\Lambda ) (\iii^{(1)}-\rr,\iii^{(2)}) =  \sum_{\iii}\PP^* (\iii^{(1)}-\rr,\iii)  \Lambda(\iii,\iii^{(2)}) 
 =  \sum_{\iii} \PP^* (\iii^{(1)}-\rr,\iii) { \mathbf{1}(\iii^{(2)}\preceq \iii)  \over \rho(\iii)}. $$
 Note that for a given $\iii^{(1)}-\rr\in\EE^*$ the only nonzero entries of $\PP^* (\iii^{(1)}-\rr,\iii)$ are for $\iii=\iii^{(1)}-\rr$ or $\iii=\iii^{(1)}-\rr+\s_j$ (if $\iii\in\EE^*$), where $\s_j=(0,\ldots,0,1,0,\ldots,0)$ (with 1 at position $j$).
 We have
 \medskip\par 
$(\PP^*\Lambda ) (\iii^{(1)}-\rr,\iii^{(2)})=$
 $$
 \begin{array}{lrll}
\displaystyle\PP^* (\iii^{(1)}-\rr,\iii^{(1)}-\rr) {\mathbf{1}(\iii^{(2)}\preceq \iii^{(1)}-\rr)\over \rho(\iii^{(1)}-\rr)}  +
 \sum_{j: i_j^{(1)}-r_j<N_j} \PP^* (\iii^{(1)}-\rr,\iii^{(1)}-\rr+\s_j) {\mathbf{1}(\iii^{(2)}\preceq \iii^{(1)}-\rr+\s_j)\over \rho(\iii^{(1)}-\rr+\s_j)}  \\[18pt]
 = \displaystyle \left(1-\sum_{j: i_j^{(1)}-r_j<N_j} p_j\right) {\mathbf{1}(\iii^{(2)}\preceq \iii^{(1)}-\rr)\over \rho(\iii^{(1)}-\rr)}  
   +\displaystyle \sum_{j: i_j^{(1)}-r_j<N_j}  {\mathbf{1}(\iii^{(2)}\preceq \iii^{(1)}-\rr+\s_j)\over \rho(\iii^{(1)}-\rr+\s_j)}p_j,  \\[18pt]
 \end{array}
 $$
%
%
\noindent 
thus
$$\PP(\iii^{(1)},\iii^{(2)}) 
=\sum_{\rr=(r_1,\ldots,r_d)\in\{0,1\}^d:\ \iii^{(1)}-\rr\in\EE^*} (-1)^{|\rr|}\rho(\iii^{(1)}-\rr) $$
\begin{equation}\label{eq:hm1}
 \times \left[\left(1-\sum_{j: i_j^{(1)}-r_j<N_j} p_j\right) {\mathbf{1}(\iii^{(2)}\preceq \iii^{(1)}-\rr)\over \rho(\iii^{(1)}-\rr)} 
+\sum_{j: i_j^{(1)}-r_j<N_j}  { \mathbf{1}(\iii^{(2)}\preceq \iii^{(1)}-\rr+\s_j)\over \rho(\iii^{(1)}-\rr+\s_j)}p_j \right].
\end{equation}
For convenience, define 
$$\begin{array}{llll}
H_1(\iii^{(1)},\iii^{(2)},\rr)&:=&\displaystyle \left(1-\sum_{j: i_j^{(1)}-r_j<N_j} p_j\right) {\mathbf{1}(\iii^{(2)}\preceq \iii^{(1)}-\rr)\over \rho(\iii^{(1)}-\rr)},\\[20pt]
H_2(\iii^{(1)},\iii^{(2)},\rr)&:=&\displaystyle \sum_{j: i_j^{(1)}-r_j<N_j}  { \mathbf{1}(\iii^{(2)}\preceq \iii^{(1)}-\rr+\s_j)\over \rho(\iii^{(1)}-\rr+\s_j)}p_j. \\   
  \end{array}
$$
\noindent
Consider  cases:

\begin{itemize}
\item \textbf{Case 1}. \textsl{Increasing some  coordinates}:  $\iii^{(2)}=\iii^{(1)}+m_1\s_{k_1}+\ldots
+m_t\s_{k_t},$ where $1\leq k_i\leq d, i=1,\ldots, t$ are $t\geq 1$ distinct integers 
and $m_i\geq 1, i=1,\ldots, t$. 
When $t\geq 2$, then indicators in both, $H_1$ and $H_2$ are equal to 0.
When $t=1$, then the indicator in $H_1$ is equal to 0, whereas 
the indicator in $H_2$ can be nonzero only in case $m=1$,
$j=k$ and $\rr=(0,\ldots,0)$. Then we have 
$$ \PP(\iii^{(1)},\iii^{(1)}+\s_k)  = (-1)^0\rho(\iii^{(1)}) {p_k\over \rho(\iii^{(1)}+\s_k)}=
{i_k^{(1)}+1\over i_k^{(1)}+2}p_k.$$

%
%
%
%
\medskip\par
\item \textbf{Case 2}. \textsl{Increasing two or more coordinates and decreasing any number of coordinates}: 
because  of the same reasons as in previous case, indicators in both, $H_1$ and $H_2$ are equal to 0.
\medskip\par
 \item \textbf{Case 3}. \textsl{Decreasing some  coordinates}:  $\iii^{(2)}=\iii^{(1)}-m_1\s_{k_1}-\ldots -m_t\s_{k_t}$,
where $1\leq k_i\leq d, i=1,\ldots,t$ are  $t\geq 1$ distinct integers and 
$1\leq m_i \leq i_{k_i},  i=1,\ldots,t$.

\par \noindent 
Let $\kkappa=(\kappa_1,\ldots,\kappa_d)$, where 
$\kappa_{k_i}=1, i=1,\ldots,t$ and $\kappa_{j}=0$ for $j\notin\{k_1,\ldots,k_t\}$.
In (\ref{eq:hm1}) we sum over all $\rr\in \{0,1\}^d$ such that $\iii^{(1)}-\rr\in\EE^*$.
Let us split this sum into two sums over disjoint sets $I_1$ and $I_2$, where 
$$\ I_1:=\{\e\in \{0,1\}^d: \e\preceq \kkappa, \iii^{(1)}-\e\in\EE^*\},\qquad
I_2:=\{\e\in  \{0,1\}^d: \e\npreceq \kkappa, \iii^{(1)}-\e\in\EE^*\}.$$

\noindent
Consider $\rr'=(r_1',\ldots,r_d')\in I_2$. Since it is incomparable with $\kkappa$ it means that 
for some $q\geq 1$ we have $a_1,\ldots,a_q$ such that $\{a_1,\ldots,a_q\}\cap \{k_1,\ldots, k_t\}=\emptyset$
and $r_{a_i}'=1, i=1,\ldots,q$. 
Then the  indicator in $H_1$ is equal to 0. 
The second indicator can be nonzero only when $q=1$ and $j=a_1$. 
Thus, for any $\rr\in I_1$ we have that $\rr + s_n \in I_2,$ for all 
$1\leq n\leq d$ such that $n\neq k_i, i=1,\ldots,t$. We  have

$$\sum_{\rr \npreceq \kkappa: \iii^{(1)}-\rr\in\EE^*}
  (-1)^{|\rr|}  \rho(\iii^{(1)}-\rr)\left(H_1(\iii^{(1)},\iii^{(2)},\rr)+H_2(\iii^{(1)},\iii^{(2)},\rr)\right)=$$
$$\sum_{n: 1\leq n\leq d\atop n\neq \kappa_i, i=1,\ldots,t}\sum_{\rr \preceq \kkappa: \iii^{(1)}-\rr\in\EE^*}
  (-1)^{|\rr+\s_n|}  \rho(\iii^{(1)}-\rr-\s_n)\left(\sum_{j: i_j^{(1)}-r_j<N_j}  { \mathbf{1}(\iii^{(1)}-m_1\s_{k_1}-\ldots -m_t\s_{k_t}\preceq \iii^{(1)}-\rr-\s_n+\s_j)\over \rho(\iii^{(1)}-\rr-\s_n+\s_j)}p_j \right).$$
  \medskip\par\noindent  The indicator is nonzero only when $j=n$, and $r_n=0$ for $n\notin\{k_1,\ldots,k_t\}$, thus
$$\begin{array}{llll}
   &=&\displaystyle\sum_{\rr \preceq \kkappa: \iii^{(1)}-\rr\in\EE^*}
  (-1)^{|\rr|+1} \sum_{n: 1\leq n\leq d\atop n\neq \kappa_i, i=1,\ldots,t}
  \left(\mathbf{1}(i_n^{(1)}-r_n<N_n) 
  {  \rho(\iii^{(1)}-\rr-\s_n)\over \rho(\iii^{(1)}-\rr)}p_j \right)\\[12pt]
  &=&\displaystyle-\sum_{\rr \preceq \kkappa: \iii^{(1)}-\rr\in\EE^*}
  (-1)^{|\rr|} \sum_{n: 1\leq n\leq d\atop n\neq \kappa_i, i=1,\ldots,t}   
  {  i_n^{(1)}\over i_n^{(1)}+1}p_j =0,
  \end{array}
$$ 
  since the second sum does not depend on $\rr$.
\medskip\par \noindent
Consider $\rr\in I_1$. Then indicators in both $H_1$ and $H_2$ are nonzero, we have 
$$\begin{array}{llll}
S_1&:=&\displaystyle\sum_{\rr \preceq \kkappa: \iii^{(1)}-\rr\in\EE^*} (-1)^{|\rr|}\rho(\iii^{(1)}-\rr)\left(H_1(\iii^{(1)},\iii^{(2)},\rr)+H_2(\iii^{(1)},\iii^{(2)},\rr)\right)    \\[14pt]
&=&\displaystyle\sum_{\rr \preceq \kkappa: \iii^{(1)}-\rr\in\EE^*}(-1)^{|\rr|}\left[ \left(1-\sum_{j: i_j^{(1)}-r_j<N_j} p_j\right) {\rho(\iii^{(1)}-\rr)\over \rho(\iii^{(1)}-\rr)} +\sum_{j: i_j^{(1)}-r_j<N_j
}  {\rho(\iii^{(1)}-\rr)\over \rho(\iii^{(1)}-\rr+\s_j)}p_j\right] \\[14pt]
&=&\displaystyle\sum_{\rr \preceq \kkappa: \iii^{(1)}-\rr\in\EE^*}(-1)^{|\rr|}\left[1-\sum_{j: i_j^{(1)}-r_j<N_j} p_j +\sum_{j: i_j^{(1)}-r_j<N_j}
\left\{ 
{i_{j}^{(1)}+1\over i_{j}^{(1)}+2} \mathbf{1}(r_j=0)
+{i_{j}^{(1)}\over i_{j}^{(1)}+1} \mathbf{1}(r_j=1)
\right\}p_j\right].
  \end{array}
$$
\noindent 
Consider cases:
\begin{itemize}
\item[$a)$] $t=1$, i.e., we decrease only one coordinate. In this case 
$\kkappa=(0,\ldots,0,1,0,\ldots,0)$ with only one 1 at position $k$.
Thus there are only two $\rr$ such that $\rr\preceq \kkappa,$ namely 
$\rr=(0,\ldots,0)$ or $\rr=\kkappa$. We have 
$$\begin{array}{llll}
S_1&=& &\displaystyle
\left[1-\sum_{j: i_j^{(1)}-0<N_j} p_j +\sum_{j: i_j^{(1)}-0<N_j}
\left\{ 
{i_{j}^{(1)}+1\over i_{j}^{(1)}+2}
\right\}p_j\right]  \\[16pt]
& &
-&
\displaystyle\left[1-\sum_{j: i_j^{(1)}-\mathbf{1}(j=k)<N_j} p_j +\sum_{j: i_j^{(1)}-\mathbf{1}(j=k)<N_j}
\left\{ 
{i_{j}^{(1)}+1\over i_{j}^{(1)}+2}\mathbf{1}(j\neq k)
+{i_{j}^{(1)}\over i_{j}^{(1)}+1}\mathbf{1}(j=k)
\right\}p_j\right].
\end{array}
$$
Note that for $j\neq k$ all the corresponding terms (for $\rr=(0,\ldots,0)$ and $\rr=\kkappa$)
are the same, thus they sum up to 0. The remaining terms:
$$\begin{array}{llll}
S_1&= &\displaystyle
\left(p_k-{i_k^{(1)}\over i_k^{(1)}+1} p_k\right)\mathbf{1}(i_k-1<N_k)
+\left(-p_k+{i^{(1)}_k +1\over i^{(1)}+2}p_k\right)\mathbf{1}(i_k<N_k)\\[16pt]
&=&\displaystyle{1\over i_k^{(1)}+1} p_k 
-{1\over i^{(1)}+2}p_k\mathbf{1}(i_k<N_k).
\end{array}
$$
\noindent
Finally, we have 
$$
\PP(\iii^{(1)},\iii^{(1)}-m\cdot \s_k)
=\left\{ \begin{array}{lllllll}
  \displaystyle {1\over (i_k^{(1)}+1) (i_k^{(1)}+2)} p_k & \textrm{if}\ i^{(1)}_k<N_k, \\[15pt]
  \displaystyle {1\over N_k+1} p_k & \textrm{if}\ i^{(1)}_k=N_k. \\
 \end{array}\right.
$$
\item[$b)$] $t\geq 0$. Things are different in this case.
Consider $\rr=(r_1,\ldots,r_d)\preceq \kkappa$ and fixed $r_n$, where  $n\in\{k_1,\ldots,k_t\}$. Then there are $2^{t-1}$ 
different $\rr$ in $S_1$, from which exactly 
$2^{t-2}$ gives $(-1)^{|\rr|}=1$ and exactly $2^{t-2}$ gives  $(-1)^{|\rr|}=1$,
resulting in vanishing the terms ${i_j^{(1)}+1\over i_j^{(1)}+2}$ or 
${i_j^{(1)}\over i_j^{(1)}+1}$ (depending on the value of $r_n$). This implies that 
$S_1=0$. For example, for $t=2$ and, for simplicity, for $d=2$, there are four following terms in $S_1$:

$$
\begin{array}{lllll}
\rr=(0,0): & \quad & & \displaystyle \mathbf{1}(i_1^{(1)}-0<N_1){i_1^{(1)}+1\over i_1^{(1)}+2} 
+ \mathbf{1}(i_1^{(1)}-0<N_1){i_2^{(1)}+1\over i_2^{(1)}+2}, \\[16pt]
\rr=(0,1): & \quad &-\bigg[ & \displaystyle \mathbf{1}(i_1^{(1)}-0<N_1){i_1^{(1)}+1\over i_1^{(1)}+2} 
+ \mathbf{1}(i_1^{(1)}-1<N_1){i_2^{(1)}\over i_2^{(1)}+1}&\bigg],\\[16pt]
\rr=(1,0): & \quad &-\bigg[ & \displaystyle \mathbf{1}(i_1^{(1)}-1<N_1){i_1^{(1)}\over i_1^{(1)}+1} 
+ \mathbf{1}(i_1^{(1)}-0<N_1){i_2^{(1)}+1\over i_2^{(1)}+2}&\bigg],\\[16pt]
\rr=(1,1): & \quad & & \displaystyle \mathbf{1}(i_1^{(1)}-1<N_1){i_1^{(1)}\over i_1^{(1)}+1} 
+ \mathbf{1}(i_1^{(1)}-1<N_1){i_2^{(1)}\over i_2^{(1)}+1},
\end{array}
$$
which sum up to 0.
\par\noindent
\textsl{Remark:} In case $t=1$ for fixed $r_{k_1}$ there was no corresponding $n\neq k_1$ which could 
make the terms vanish.

\end{itemize}
\smallskip\par

\item  \textbf{Case 4}. \textsl{Increasing one, decreasing another coordinate:}  $\iii^{(2)}=\iii^{(1)} -m_1 \s_{k_1}+m_2\s_{k_2}$. 
We have shown that increasing/decreasing $t\geq 2$ coordinates has probability 0,
thus there is no need to consider the case where we increase and decrease any number of coordinates 
in one step. 
\par \noindent 
In this case the indicator in $H_1$ is zero. Concerning $H_2$. Let,
$\kkappa=(0,\ldots,0,1,0,\ldots,0)$ with one 1 at position $k_1$. Note that for $\rr\npreceq \kkappa$,
the indicator in $H_2$ is also 0. Thus, the only nonzero terms are for either 
$\rr=(0,\ldots,0)$ or $\rr=\kkappa$ (and then $j=k_2$):
$$
\begin{array}{lllll}
\rr=(0,\ldots,0): & \quad  & \displaystyle \mathbf{1}(i_{k_1}^{(1)}-0<N_1){\rho(\iii_1^{(1)})\over \rho(\iii_1^{(1)}+\s_{k_2})}, \\[16pt]
\rr=\kkappa: & \quad & - \displaystyle \mathbf{1}(i_{k_1}^{(1)}-0<N_1){\rho(\iii_1^{(1)})\over \rho(\iii_1^{(1)}+\s_{k_2})}, \\[16pt]
\end{array}
$$
what sums up to 0.
\smallskip\par
\item \textbf{Case 5}.  \textsl{Staying at the same state:}  $\iii^{(2)}=\iii^{(1)}$. Then 
the indicator 
$\mathbf{1}(\iii^{(2)}\preceq \iii^{(1)}-\rr)$ is nonzero only when $\rr=(0,\ldots,0)$,
whereas the indicator $ \mathbf{1}(\iii^{(2)}\preceq \iii^{(1)}-\rr+\s_j) $
is nonzero when $\rr=(0,\ldots,0)$ and any $j=1,\ldots,d$ or when $\rr=\s_j$. 
We have 
$$
\begin{array}{lllllllllll} 
\PP(\iii^{(1)},\iii^{(1)})& =& \displaystyle (-1)^0 \rho(\iii^{(1)})    \displaystyle \left[\left(1-\sum_{j: i^{(1)}_j-0<N_j} p_j\right) {1\over \rho(\iii^{(1)})} + \sum_{j: i^{(1)}_j-0<N_j}  {p_j\over \rho(\iii^{(1)}+\s_j)} \right]\\[18pt]
    & & \displaystyle -\sum_{k: i^{(1)}_k-1\geq 0}\rho(\iii^{(1)}-\s_k)  \displaystyle\left[ 
 \sum_{j: i^{(1)}_j-\mathbf{1}(j=k)<N_j}  {p_j\over \rho(\iii^{(1)}-\s_k+\s_j)} \mathbf{1}(\iii^{(1)}\preceq \iii^{(1)}-\s_k+\s_j)
 \right]\\[18pt]

 & = &\displaystyle 1-\sum_{j: i^{(1)}_j<N_j} p_j  + \sum_{j: i^{(1)}_j<N_j}  {\rho(\iii^{(1)})\over \rho(\iii^{(1)})+\s_j)}p_j  
 - \sum_{k: i^{(1)}_k\geq 1}    {\rho(\iii^{(1)}-\s_k)\over \rho(\iii^{(1)})}p_k \\[24pt] 
&=& \displaystyle 1-\sum_{j: i^{(1)}_j<N_j} \left(1-{i^{(1)}_j+1\over i^{(1)}_j+2} \right)p_j
- \sum_{k=1}^d  {i^{(1)}_k\over i^{(1)}_k+1} p_k \\[18pt]
&=& \displaystyle 1-\sum_{j: i^{(1)}_j<N_j} \left({1\over i^{(1)}_j+2}\right)p_j
- \sum_{j: i^{(1)}_j<N_j }  {i^{(1)}_j\over i^{(1)}_j+1} p_j- \sum_{j:  i^{(1)}_j=N_j}  {i^{(1)}_j\over i^{(1)}_j+1} p_j\\[18pt]
&=& \displaystyle 1-\sum_{j: i^{(1)}_j<N_j} \left(1-{1\over (i^{(1)}_j+1)(i^{(1)}_j+2)}\right)p_j
- \sum_{j: i^{(1)}_j=N_j}  {N_j\over N_j+1} p_j.
\end{array}
$$
\end{itemize}
The assumption (\ref{eq:lattice_assumptions}) implies that $\PP(\iii^{(1)},\iii^{(1)})\geq 0$.
We have considered all the transitions.
Let us check that each row of calculated $\PP$ sums up to 1.
We have (with the convention $\sum_{m=1}^0 f(m)\equiv 0$)
\bigskip\par 
$\displaystyle\sum_{\iii^{(2)}\in\EE^*}\PP(\iii^{(1)},\iii^{(2)})=$
\vspace{-0.75cm}
$$
\begin{array}{lllllllllll} 
 &  &\displaystyle \PP(\iii^{(1)},\iii^{(1)})+\sum_{j:  i^{(1)}_j<N_j }\PP(\iii^{(1)},\iii^{(1)}+\s_j)
+\sum_{j=1}^d \sum_{m=1}^{i_j^{(1)}}\PP((\iii^{(1)},\iii^{(1)}-m\cdot\s_j)\\[18pt] 
 & = &\displaystyle 1-\sum_{j: i^{(1)}_j<N_j} \left(1-{1\over (i^{(1)}_j+1)(i^{(1)}_j+2)}\right)p_j
- \sum_{j: i^{(1)}_j=N_j}  {N_j\over N_j+1} p_j\\[18pt] 
 &  &\displaystyle +\sum_{j: i^{(1)}_j<N_j} { i^{(1)}+1\over  i^{(1)}+2}p_j+\sum_{j: i^{(1)}_j<N_j}  \sum_{m=1}^{i_j^{(1)}}\PP((\iii^{(1)},\iii^{(1)}-m\cdot\s_j)
+\sum_{j: i^{(1)}_j=N_j}  \sum_{m=1}^{i_j^{(1)}}\PP((\iii^{(1)},\iii^{(1)}-m\cdot\s_j)\\[18pt] 
 & = &\displaystyle 1-\sum_{j: i^{(1)}_j<N_j} \left(1-{1\over (i^{(1)}_j+1)(i^{(1)}_j+2)} - { i^{(1)}_j+1\over  i^{(1)}_j+2}\right)p_j
- \sum_{j: i^{(1)}_j=N_j}  {N_j\over N_j+1} p_j\\[18pt] 
 &  &\displaystyle +\sum_{j: i^{(1)}_j<N_j}  {i_j^{(1)}\over (i_j^{(1)}+1)(i_j^{(1)}+2) }
+\sum_{j: i^{(1)}_j=N_j}  {N_j\over N_j+1} p_j\\[18pt] 
 & = &\displaystyle 1-\sum_{j: i^{(1)}_j<N_j} \left(1-{1\over (i^{(1)}_j+1)(i^{(1)}_j+2)} - { i^{(1)}_j+1\over  i^{(1)}_j+2} -{i_j^{(1)}\over (i_j^{(1)}+1)(i_j^{(1)}+2) } \right)p_j=1
\end{array}
$$
 
\end{proof}


%
%
%

\begin{proof}[Proof of Theorem \ref{twr:coupon_antidual_ver2}]
Note that $(0,\ldots,0)$ is the minimal state, and $X^*$ starts at this state $\nu^*=\delta_{(0,\ldots,0)}$, thus - by Remark \ref{rem:star_min} - this 
is also the initial distribution of the antidual chain, i.e., $\nu=\nu^*$.\par
\noindent
For convenience, define
$$f(\iii,k)={\sum_{\iii' \preceq \iii} \pi(\iii')\over \sum_{\iii''\preceq \iii+\s_k}\pi(\iii'')} \quad \mathrm{for}\quad \iii: i_k=0.$$
\noindent
For the stationary distribution $\pi$ given in (\ref{eq:gen_coupon_coll_antidual_pi}) we have
$${\pi(\iii^{(1)}+\s_k)\over \pi(\iii^{(1)})} = {a_k\over 1-a_k}, \quad {\pi(\iii^{(1)}-\s_k)\over \pi(\iii^{(1)})} = {1-a_k\over a_k},$$
$$f(\iii,k) = {\displaystyle  \sum_{\iii'\preceq \iii} \prod_{j=1}^d[a_j \mathbf{1}(i_j'=1)+(1-a_j) \mathbf{1}(i'_j=0)]\over
\displaystyle  \sum_{\iii''\preceq \iii+\s_k} \prod_{j=1}^d[a_j \mathbf{1}(i''_j=1)+(1-a_j) \mathbf{1}(i''_j=0)]}.
$$
Denote
$$\xi(\iii,k)=\prod_{j=1\atop j\neq k}^d  [a_j \mathbf{1}(i_j=1)+(1-a_j) \mathbf{1}(i_j=0)].$$
The sum in denominator of $f(\iii,k)$ can be split into two sums: for $\iii'': i''_k=0$ and $\iii'':i''_k=1$. We have 
$$ f(\iii,k)= {\displaystyle \sum_{\iii'\preceq \iii} \xi(\iii',k) (1-a_k)\over
\displaystyle \sum_{\iii''\preceq \iii+\s_k \atop \iii''_k=0} \xi(\iii'',k)(1-a_k) + \sum_{\iii''\preceq \iii+\s_k \atop i''_k=1}\xi(\iii'',k) a_k}=1-a_k.$$

\noindent
Let us proceed with $\widehat{\PP^*}$.
$$\widehat{\PP^*}(\iii^{(2)},\iii^{(1)})={(\boldsymbol{\pi}\C)(\iii^{(2)})\over(\boldsymbol{\pi}\C)(\iii^{(1)})}\PP^*(\iii^{(2)},\iii^{(1)})=
\left\{
\begin{array}{llllll}
  f(\iii^{(2)},k) p_k = (1-a_k)p_k& \mathrm{if} & \iii^{(1)}=\iii^{(2)}+\s_k,\\[8pt]
 \displaystyle 1-\sum_{j: i^{(2)}_j=0} p_j & \mathrm{if} & \iii^{(1)}=\iii^{(2)}.\\
\end{array}\right.
$$
Note that $\widehat{\PP^*}$ is not a stochastic matrix, since we have 
\begin{equation*}\label{eq:sumf}
\sum_\iii \widehat{\PP^*}(\iii^{(2)},\iii)=\sum_{j:i^{(2)}_j=0} f(\iii^{(2)},j)p_j + 1-\sum_{j:i^{(2)}_j=0} p_j=1-\sum_{j:i^{(2)}_j=0}\left(1-f(\iii^{(2)},j)\right)p_j
=1-\sum_{j:i^{(2)}_j=0}a_j p_j<1.
\end{equation*}
Now, calculating the antidual chain from Theorem \ref{twr:main}, we have 
$$\PP(\iii^{(1)},\iii^{(2)})={\pi(\iii^{(2)})\over\pi(\iii^{(1)})} ( (\C^T)^{-1} \widehat{\PP^*} \C^T)(\iii^{(1)},\iii^{(2)}) = {\pi(\iii^{(2)})\over \pi(\iii^{(1)})}(\C(\widehat{\PP^*})^T\C^{-1})(\iii^{(2)},\iii^{(1)}) $$
\begin{equation}\label{eq:antidual_cube_tr}
= {\pi(\iii^{(2)})\over\pi(\iii^{(1)})}  \sum_{\iii\preceq\iii^{(1)}} \widehat{\PP^*}(\iii,\{\iii^{(2)}\}^\uparrow)(-1)^{|\iii^{(1)}-\iii|}, 
\end{equation}
where we applied the   M\"obius function  for this ordering: $\C^{-1}(\iii,\iii^{(1)})= (-1)^{|\iii^{(1)}-\iii|} \mathbf{1}(\iii\preceq \iii^{(1)})$
(a consequence of (\ref{eq:coord_wise_mobius_fun})).
We proceed with (\ref{eq:antidual_cube_tr}) by considering cases:
\begin{itemize}
 \item \textbf{Case 1.} \textsl{Increasing some  coordinates}: 
 $\iii^{(2)}=\iii^{(1)}+\s_{k_1}+\ldots \s_{k_t} $ for some distinct $t\geq 1$ integers 
 $1\leq k_i\leq d, i=1,\ldots,d$.
 \smallskip\par\noindent
 First note that if $t\geq 2$, than, for any $\iii\preceq\iii^{(1)}$ we
 have  $\widehat{\PP^*}(\iii,\{\iii^{(1)}+\s_{k_1}+\ldots \s_{k_t}\}^\uparrow)=0$, thus $\PP(\iii^{(1)},\iii^{(1)}+\s_{k_1}+\ldots \s_{k_M})=0$.
 \smallskip\par\noindent
 For $t=1$  the sum in (\ref{eq:antidual_cube_tr}) is following $ \sum_{\iii\preceq\iii^{(1)}} \widehat{\PP^*}(\iii,\{\iii^{(1)}+\s_k\}^\uparrow)(-1)^{|\iii^{(1)}-\iii|}$,
 the only nonzero term is for $\iii=\iii^{(1)}$, thus 
 \begin{equation*}\label{eq:case1}
 \begin{array}{llll}
\PP(\iii^{(1)},\iii^{(1)}+\s_k)&=&\displaystyle{\pi(\iii^{(1)}+\s_k)\over\pi(\iii^{(1)})} \widehat{\PP^*}(\iii^{(1)},\{\iii^{(1)}+\s_k\}^\uparrow)={\pi(\iii^{(1)}+\s_k)\over\pi(\iii^{(1)})}  f(\iii^{(1)},k)p_k\\[14pt]
  &=&\displaystyle{a_k\over 1-a_k} (1-a_k)p_k=a_kp_k.   
 \end{array}
 \end{equation*}

 \medskip\par 
 \item \textbf{Case 2.} \textsl{Increasing two or more coordinates and decreasing any number of coordinates}: 
 because of the same reasons as in previous case (we would have to increase at 
 least two coordinates in one step) such transition has probability 0.
 \medskip\par

 \item \textbf{Case 3:}
 $\iii^{(2)}=\iii^{(1)}-\s_{k_1}-\ldots-\s_{k_t}, t\geq 1$.
 Let us split $\{\e\preceq \iii^{(1)}\}$ into five disjoint sets:
 
 $$I_1=\{\iii^{(1)}\}, \quad I_2=\{\iii^{(2)}\},\quad  I_3=\{\e: \e\prec \iii^{(2)}\}, \quad
 I_4=\{\e: \iii^{(2)}\prec \e \prec \iii^{(1)}\}, \quad  I_5=\{\e: \e\npreceq \iii^{(2)}\},
 $$
 where $\e\prec\e'$ means that $\e\preceq \e'$ and $\e\neq \e'$, and 
 $\e\npreceq \e'$ means that $\e$ and $\e'$ are incomparable.
 Define also
 $$S_m:=\sum_{\iii\in I_m} \widehat{\PP^*}(\iii,\{\iii^{(2)}\}^\uparrow)(-1)^{|\iii^{(1)}-\iii|}, \quad m=1,2,3,4,5.$$
 We have
  $$ 
 \begin{array}{llll}
  S_1
  &=&\displaystyle \widehat{\PP^*}(\iii^{(1)},\{\iii^{(2)}\}^\uparrow)=1-\sum_{j:i_j^{(1)}=0}a_jp_j.\\[18pt]
  
  \displaystyle S_2
  &=&\displaystyle \widehat{\PP^*}(\iii^{(2)},\{\iii^{(2)}\}^\uparrow)(-1)^{|\iii^{(1)}-\iii^{(2)}|}=\left(1-\sum_{j:i_j^{(1)}=0}a_jp_j-\sum_{j\in\{k_1,\ldots,k_t\}}a_jp_j\right)(-1)^{|\iii^{(1)}-\iii^{(2)}|}.\\[18pt]
  
  \displaystyle S_3
  &=&\displaystyle 
  \sum_{j:i^{(1)}_j=1\atop j\notin\{k_1,\ldots,k_t\} }
  \widehat{\PP^*}(\iii^{(2)}-\s_j,\{\iii^{(2)}\}^\uparrow)(-1)^{|\iii^{(1)}-\iii|}=
  (-1)^{|\iii^{(1)}-\iii^{(2)}-1|}\sum_{j:i^{(1)}_j=1\atop j\notin\{k_1,\ldots,k_t\} } (1-a_j)p_j.\\[20pt]
  
  \displaystyle S_4 
  &=&\displaystyle 
  \sum_{\iii\in I_4}\left(1-\sum_{j:i_j=0} a_jp_j\right)(-1)^{|\iii^{(1)}-\iii|}=\sum_{\iii\in I_4}\left(
  1-\sum_{j:i^{(1)}_j=0} a_jp_j - 
 \sum_{j\in\{k_1,\ldots,k_j\}\atop  i_j=0} a_jp_j\right)(-1)^{|\iii^{(1)}-\iii|}.\\[20pt]

  \displaystyle S_5
  
  &=&\displaystyle 
  \sum_{\emptyset\neq\{b_1,\ldots,b_z\}\subseteq\{k_1,\ldots,k_t\}}
  \sum_{j: i^{(1)}_j=1\atop j\notin\{k_1,\ldots,k_t\}} 
  
  \widehat{\PP^*}(\iii^{(1)}-\s_{b_1}-\ldots,\s_{b_z}-\s_j,\{\iii^{(2)}\}^\uparrow)(-1)^{t-z+1}
  \\[22pt]

  &=&\displaystyle 
  \sum_{\emptyset\neq\{b_1,\ldots,b_z\}\subseteq\{k_1,\ldots,k_t\}}
  \sum_{j: i^{(1)}_j=1\atop j\notin\{k_1,\ldots,k_t\}}  (1-a_j)p_j(-1)^{t-z+1}
  =(-1)^t\sum_{j: i^{(1)}_j=1\atop j\notin\{k_1,\ldots,k_t\}}  (1-a_j)p_j.\\  
  \end{array}
$$
\noindent
Let us consider cases $t=1$ and $t\geq 2$ separately.

\begin{itemize}
\item[$a)$] $t=1$, i.e., $\iii^{(2)}=\iii^{(1)}-\s_{k}\  (k_1\equiv k)$.
Note that then  $I_4 = \emptyset $. We have 

  $$ 
 \begin{array}{llll}
  \displaystyle S_1
  &=&\displaystyle 1-\sum_{j:i_j^{(1)}=0}a_jp_j.\\[18pt]
  
  \displaystyle S_2
  &=&\displaystyle -\left(1-\sum_{j:i_j^{(1)}=0}a_jp_j-a_{k}p_{k}\right).\\[18pt]
  
  \displaystyle S_3
  &=&\displaystyle 
 \sum_{j:i^{(1)}_j=1\atop j\neq k } (1-a_j)p_j.\\[20pt]
     
  \displaystyle S_5
  &=&\displaystyle 
  -\sum_{j:i_j^{(1)}=1, j\neq k} \widehat{\PP^*}(\iii^{(1)}-\s_j,\{\iii^{(2)}\}^\uparrow)=- \sum_{j:i_j^{(1)}=1, j\neq k} ( 1-a_j)p_j. \\[14pt]
  \end{array}
$$
We have $S_1+S_2 +S_3 + S_4 +S_5 = a_kp_k$ and finally
$$\PP(\iii^{(1)},\iii^{(1)}-\s_k)={\pi(\iii^{(1)}-\s_k)\over \pi(\iii^{(1)})} (S_1+S_2+S_3 +S_4)={1-a_k\over a_k} a_kp_k =(1-a_k)p_k.$$

\item[$b)$] $t\geq 2$. Consider first $t=2$.
Assume thus that  $\iii^{(2)}=\iii^{(1)}-\s_{k_1}-\s_{k_2}.$
We have 

  $$ 
 \begin{array}{llll}
  S_1
  &=&\displaystyle 1-\sum_{j:i_j^{(1)}=0}a_jp_j.\\[18pt]
  
  \displaystyle S_2
  &=&\displaystyle 
  1-\sum_{j:i_j^{(1)}=0}a_jp_j-a_{k_1}p_{k_1}-a_{k_2}p_{k_2}.\\[18pt]
  
  \displaystyle S_3
  &=&\displaystyle 
  -\sum_{j:i^{(1)}_j=1\atop j\notin\{k_1,k_2\} } (1-a_j)p_j.\\[22pt]
  
  \displaystyle S_4 
  &=&\displaystyle 
 - \widehat{\PP^*}(\iii^{(1)}-\s_{k_1},\{\iii^{(1)}-\s_{k_1}-\s_{k_2}\}^\uparrow)
  -\widehat{\PP^*}(\iii^{(1)}-\s_{k_2},\{\iii^{(1)}-\s_{k_1}-\s_{k_2}\}^\uparrow)
  \\[20pt]
  
  &=&\displaystyle 
  -(1-\sum_{j:i^{(1)}_j=0} a_jp_j - 
  a_{k_2} p_{k_2}) - (1-\sum_{j:i^{(1)}_j=0} a_jp_j - 
  a_{k_1} p_{k_1}).
  \\[20pt]
  
  \displaystyle S_5
  
  &=&\displaystyle \sum_{j: i^{(1)}_j=1\atop j\notin\{k_1,\ldots,k_t\}}  (1-a_j)p_j.\\  
  \end{array}
$$
Summing up,  $S_1+S_2+S_3+S_4+S_5=0$, what is also the case for $t>2$
(the proof, although longer, is quite similar, we skip the details). This means that for $t\geq 2$
$$\PP(\iii^{(1)},\iii^{(1)}-\s_{k_1}-\ldots-\s_{k_t})=0.$$

\end{itemize}

  \item \textbf{Case 4.} \textsl{Increasing one, decreasing another coordinate:}  $\iii^{(2)}=\iii^{(1)}+\s_{k_1}-\s_{k_2}$. 
  We have shown that increasing/decreasing $t\geq 2$ coordinate has probability 0, thus it suffices to consider only changing two coordinates (one increasing,
  the other decreasing).
  Then the  the summands    $ \sum_{\iii\preceq\iii^{(1)}} \widehat{\PP^*}(\iii,\{\iii^{(1)}+\s_{k_1}-\s_{k_2}\}^\uparrow)(-1)^{|\iii^{(1)}-\iii|}$
  are nonzero only for  $\iii=\iii^{(1)}$ or $\iii=\iii^{(1)}-\s_{k_2}$, we have 
  $$
 \begin{array}{llll}
\widehat{\PP^*}(\iii^{(1)},\{\iii^{(1)}+\s_{k_1}-\s_{k_2}\}^\uparrow)(-1)^{|\iii^{(1)}-\iii^{(1)}|}
 &=&\displaystyle f(\iii^{(1)},k_{1})p_{k_{1}}=(1-a_{k_{1}})p_{k_{1}}, \\[14pt]
 \widehat{\PP^*}(\iii^{(1)}-\s_{k_2},\{\iii^{(1)}+\s_{k_1}-\s_{k_2}\}^\uparrow)(-1)^{|\iii^{(1)}-\iii^{(1)}-1|}
 &=&\displaystyle -f(\iii^{(1)}-\s_{k_2},{k_{1}})p_{k_{1}}=-(1-a_{k_{1}})p_{k_{1}}, 
 \end{array}
 $$
 thus $\PP(\iii^{(1)},\iii^{(1)}+\s_{k_1}-\s_{k_2})=0$.
 
 \smallskip\par
 \item \textbf{Case 5.} \textsl{Staying at the same state:} $\iii^{(2)}=\iii^{(1)}$. Then we have 
 $$\begin{array}{llll}
 \PP(\iii^{(1)},\iii^{(1)}) & =&\displaystyle\sum_{\iii\preceq\iii^{(1)}} \widehat{\PP^*}(\iii,\{\iii^{(1)}\}^\uparrow)(-1)^{|\iii^{(1)}-\iii|} \\[14pt]
 &=&  \displaystyle\widehat{\PP^*}(\iii^{(1)},\{\iii^{(1)}\}^\uparrow)-\sum_{j:\iii^{(1)}_j=1}\widehat{\PP^*}(\iii^{(1)}-\s_j,\{\iii^{(1)}\}^\uparrow).   
   \end{array}
$$
 First term is equal to $\sum_\iii  \widehat{\PP^*}(\iii^{(1)},\iii)$, in the latter, the only possibility is to change $j$-th 
 coordinate of $\iii^{(1)}-\s_j$ to one:
 \begin{equation*}\label{eq:case3}
 \begin{array}{lllll}
 \PP(\iii^{(1)},\iii^{(1)})&=&\displaystyle 1-\sum_{j:i^{(1)}_j=0} p_j(1-f(\iii^{(1)},j))-\sum_{j:i^{(1)}_j=1} f(\iii^{(1)}-\s_j,j)p_j\\[16pt]
  & =& \displaystyle 1-  \sum_{j:i^{(1)}_j=0} a_j p_j -  \sum_{j:i^{(1)}_j = 1} (1-a_j)p_j.   \\
  \end{array}
 \end{equation*}

\end{itemize}
\medskip\par 
Finally, we obtain matrix $\PP$ given in (\ref{eq:antidual_N1}).
\end{proof}
\begin{remark}\label{rm:two_coord}\rm
 Showing that $\PP(\iii^{(1)},\iii^{(1)}+\s_{k_1}-\s_{k_2})=0$
  relied heavily on the fact that for the stationary distribution given in (\ref{eq:gen_coupon_coll_antidual_pi}),
 we had $f(\iii,j)=1-a_j$ and it did not depend on $\iii$. That is why the terms 
 $f(\iii^{(1)},{k_{1}})p_{k_{1}}$ and $f(\iii^{(1)},{k_{1}})p_{k_{1}}$
cancelled out. Similarly, it is the reason why decreasing $t\geq 2$
 coordinates has probability 0.
For other, not product-form stationary distributions,
such transitions are possible.
\end{remark}

\subsection{Proof of Theorem \ref{thm:FSST_gen}}\label{sec:proof_FSST}
\noindent Let $\X^*$ be an absorbing chain on $\EE=\{1,\ldots,M\}, M\geq 2$ with transition matrix:
$$ \PP^*(k,s) 
= \left\{ 
\begin{array}{lll}
 p_k & \textrm{if\ } s=k+1, \\[10pt]
 1-p_k & \textrm{if\ } s=k, 
\end{array}\right. 
 $$
where, for convenience, we set $p_M=0$. Let $\nu^*=(a_1,\ldots,a_M)$ be its initial distribution. 
This is a pure birth chain, thus its absorption time $T^*$ is distributed as (\ref{eq:FSST_T}).
We will show that $\PP$ is its sharp antidual chain. \par \noindent
We consider the total ordering $\preceq:=\leq$. Then the link given in (\ref{eq:link0}) reads
$$\Lambda(k,s)={\pi(s)\mathbf{1}(s\leq k)\over H(k)}.$$
The inverse  $\Lambda^{-1}$ can be easily derived:
$$\Lambda^{-1}(k,s) 
=\left\{ 
\begin{array}{lll}
 \displaystyle {H(k) \over \pi(k)} & \textrm{if} & s=k, \\[12pt]
\displaystyle  -{H(k-1) \over \pi(k)} & \textrm{if} & s=k-1.
\end{array}
\right. 
$$
Let us calculate 
$$
\begin{array}{llllll}
\PP^*\Lambda(k,s)&=&\displaystyle \sum_{r}\PP^*(k,r)\Lambda(r,s)\\[10pt]
 & = & \displaystyle \PP^*(k,k)\Lambda(k,s) + \mathbf{1}(k<M)\PP^*(k,k+1)\Lambda(k+1,s) \\[10pt]
 & = & \displaystyle {\pi(s)\over H(k)}(1-p_k)\mathbf{1}(s\leq k)+{\pi(s)\over H(k+1)} p_k  \mathbf{1}(k<M) \mathbf{1}(s\leq k+1).\\[10pt]
\end{array}
$$
Calculating transitions of the antidual chain:
$$
\begin{array}{llllll}
\PP(k,s) & = & \displaystyle\Lambda^{-1}\PP^*\Lambda(k,s)=\sum_{r}\Lambda^{-1}(k,r)\PP^*\Lambda(r,s)\\[10pt]
 & = & \displaystyle{H(k) \over \pi(k)} \PP^*\Lambda(k,s) - \mathbf{1}(k>1) {H(k-1)\over \pi(k)}\PP^*\Lambda(k-1,s).
\end{array}
$$
Consider separately the cases:
 \begin{itemize}
  \item $k=1$. Then $\PP(1,s)={H(1)\over \pi(1)} \PP^*\Lambda(1,s)= \PP^*\Lambda(1,s)$.
  This is nonzero only if $s=1$ or $s=2$. 
  $$
\begin{array}{llllll}
\PP(1,1) & = & \displaystyle (1-p_1){\pi(1)\over H(1)}+p_1{\pi(1)\over H(2)}=1-p_1+p_1{\pi(1)\over \pi(1)+\pi(2)} = 1-{\pi(2)\over\pi(1)+\pi(2)}p_1, \\[10pt]
\PP(1,2) & = & \displaystyle p_1{\pi(2)\over H(2)}={ \pi(2)\over \pi(1)+\pi(2)} p_1.
\end{array}
$$

  \item $k=M$.   We have 
  $$
\begin{array}{llllll}
  \PP^*\Lambda(M,s)&=&\displaystyle(1-p_M){\pi(s)\over H(M)}=\pi(s)\\[10pt]
  \PP^*\Lambda(M-1,s)&=&\displaystyle(1-p_{M-1}) {\pi(s)\mathbf{1}(s\leq M-1)\over H(M-1)} + \mathbf{1}(M-1<M) {\pi(s)\mathbf{1}(s\leq M)\over H(M)}p_{M-1}\\[10pt]
    & = & \displaystyle(1-p_{M-1}) {\pi(s)\mathbf{1}(s\leq M-1)\over H(M-1)} +  \pi(s)p_{M-1}.
\end{array}
$$
%
%
Thus,
  $$
\begin{array}{llllll}
  \PP(M,s) & = & \displaystyle{H(M)\over \pi(M)}\PP^*\Lambda(M,s)-\mathbf{1}(M>1){H(M-1)\over \pi(M)}\PP^*\Lambda(M-1,s)\\[10pt]
   & = & \displaystyle{H(M)\over \pi(M)} \pi(s) -{H(M-1)\over \pi(M)} \left((1-p_{M-1}) {\pi(s)\mathbf{1}(s\leq M-1)\over H(M-1)} +  \pi(s)p_{M-1}\right)\\[10pt]
  & = & \displaystyle{\pi(s)\over \pi(M)} - {\pi(s)\over \pi(M)}(1-p_{M-1}) {\mathbf{1}(s\leq M-1)} -{\pi(s)\over \pi(M)}H(M-1)p_{M-1}\\[10pt]
  & = & \displaystyle {\pi(s)\over \pi(M)}\left[1 - p_{M-1} - (1-p_{M-1}) {\mathbf{1}(s\leq M-1)} +\pi(M)p_{M-1}\right]\\[10pt]
  & = & \displaystyle \left\{\begin{array}{llll}
        p_{M-1}\pi(s) & \textrm{ if } & s\leq M-1, \\
        1-p_{M-1}+p_{M-1}\pi(M)  & \textrm{ if } & s=M.
       \end{array}\right.
\end{array}
$$
%
%
%
%
%
  
  \item $1<k<M$.  We have 
  $$\PP^*\Lambda(k-1,s)=(1-p_{k-1}) {\pi(s)\mathbf{1}(s\leq k-1)\over H(k-1)} + p_{k-1} {\pi(s)\mathbf{1}(s\leq k)\over H(k)}.$$
  Thus,
    $$
\begin{array}{llllll}
  \PP(k,s)& = & \displaystyle {H(k)\over \pi(k)}\PP^*\Lambda(k,s)-{H(k-1)\over \pi(k)}\PP^*\Lambda(k-1,s)\\[12pt]
  & =& \displaystyle {H(k)\over \pi(k)}\left[(1-p_k) {\pi(s)\mathbf{1}(s\leq k)\over H(k)} +  p_k {\pi(s)\mathbf{1}(s\leq k+1)\over H(k+1)})\right]\\[12pt]
  &- & \displaystyle {H(k-1)\over \pi(k)}\left[(1-p_{k-1}) {\pi(s)\mathbf{1}(s\leq k-1)\over H(k-1)} + p_{k-1} {\pi(s)\mathbf{1}(s\leq k)\over H(k)}\right].\\[10pt]
\end{array}
$$ 
Consider three  sub-cases:  
  \begin{itemize}
  \item[$\diamond$] $s=k+1$. Then we have 
   $$\PP(k,k+1)=p_k {H(k)\over  H(k+1)}   {\pi(k+1)\over \pi(k)}. $$
  \item[$\diamond$] $s=k$. Then we have 
      $$
\begin{array}{llllll}
  \PP(k,k) & = &\displaystyle {H(k)\over \pi(k)}\left[(1-p_k) {\pi(k)\over H(k)} +  p_k {\pi(k)\over H(k+1)}\right]-{H(k-1)\over \pi(k)}\left[ p_{k-1} {\pi(k)\over H(k)}\right]\\[12pt]
  & = & \displaystyle 1-p_k   +  p_k {H(k)\over H(k+1)}-p_{k-1} {H(k-1)\over H(k)}=1-p_k\left(1 - {H(k)\over H(k+1)}\right)-p_{k-1} {H(k-1)\over H(k)}.
\end{array}
$$ 

  \item[$\diamond$] $s<k$. Then we have
     $$
     \begin{array}{llllll}
 \PP(k,s)&=& \displaystyle {H(k)\over \pi(k)}\left[(1-p_k) {\pi(s) \over H(k)} +  p_k {\pi(s) \over H(k+1)}\right]\\[12pt]
  &-&\displaystyle {H(k-1)\over \pi(k)}\left[(1-p_{k-1}) {\pi(s) \over H(k-1)} + p_{k-1} {\pi(s) \over H(k)}\right]\\[12pt]
     &=& \displaystyle  (1-p_k) {\pi(s)\over \pi(k)} +  p_k {\pi(s) \over \pi(k)} {H(k)\over H(k+1)}-(1-p_{k-1}) {\pi(s) \over \pi(k)} - p_{k-1} {\pi(s) \over \pi(k)} {H(k-1)\over H(k)}\\[12pt]
  & = & \displaystyle {\pi(s)\over \pi(k)}\left[ p_{k-1}\left( 1-{H(k-1)\over H(k)}\right)-p_{k} \left(1-{H(k)\over H(k+1)}\right) \right].
  \end{array}
$$ 
%
%
%
%

\end{itemize}

\end{itemize}
For $k\in\{1,M\}$ we obviously have $\sum_{s=1}^M\PP(k,s)=1$. For $1<k<M$ we have

   $$
     \begin{array}{llllll}
 \sum_{s=1}^M \PP(k,s)&=& \displaystyle \sum_{s=1}^{k-1}\PP(k,s) + \PP(k,k)+ \PP(k,k+1) \\[12pt]
 
  & = & \displaystyle \sum_{s=1}^{k-1} {\pi(s)\over \pi(k)}\left[ p_{k-1}\left( 1-{H(k-1)\over H(k)}\right)-p_{k} \left(1-{H(k)\over H(k+1)}\right) \right] \\[12pt]
  &  +& \displaystyle   1-p_k\left(1 - {H(k)\over H(k+1)}\right)-p_{k-1} {H(k-1)\over H(k)}
   +  p_k {H(k)\over  H(k+1)}   {\pi(k+1)\over \pi(k)} \\[12pt]
   & = & \displaystyle  {H(k-1)\over \pi(k)}\left[ p_{k-1}\left( 1-{H(k-1)\over H(k)}\right)-p_{k} \left(1-{H(k)\over H(k+1)}\right) \right] \\[12pt]
  &  +& \displaystyle   1 + p_k \left[{H(k)\over  H(k+1)}   {\pi(k+1)\over \pi(k)}-{H(k+1)-H(k)\over H(k+1)}\right] -p_{k-1} {H(k-1)\over H(k)}
    \\[12pt]
       & = & \displaystyle 1 + p_k \left[{H(k)\over  H(k+1)}   {\pi(k+1)\over \pi(k)} - {\pi(k+1)\over H(k+1)} -{H(k-1)\over\pi(k)}{\pi(k+1)\over H(k+1)}\right]\\[12pt]
       & + & \displaystyle p_{k-1} \left[{H(k-1)\over \pi(k)}\left(1-{H(k-1)\over H(k)}\right) - {H(k-1)\over H(k)}\right]\\[12pt]
       & = & \displaystyle 1 + p_k {\pi(k+1)\over H(k+1)} \left[{H(k)\over  \pi(k)} - 1  -{H(k-1)\over\pi(k)}\right] + p_{k-1} \left[{H(k-1)\over \pi(k)}-{H(k-1)\over H(k)} \left({H(k-1)\over \pi(k)}+1\right) \right]\\[12pt]
       & = & \displaystyle 1 + p_k {\pi(k+1)\over H(k+1)} \left[{H(k)-\pi(k)-H(k-1)\over  \pi(k)} \right] + p_{k-1} \left[{H(k-1)\over \pi(k)}-{H(k-1)\over H(k)} {H(k)\over \pi(k)} \right]=1.\\[12pt]
  \end{array}
$$ 
Thus (cf. (\ref{eq:sumup1}))   we considered all the cases. The only thing left to calculate is the initial distribution of the antidual chain. 
Using relation (\ref{eq:duality}) we have 
$$\nu(k)=\sum_{i=1}^M \nu^*(i)\Lambda(i,k)=\pi(k)\sum_{i=1}^M {a_i\mathbf{1}(k\leq i)\over H(i)}
=\pi(k)\sum_{i=k}^M {a_i \over H(i)}.$$
The matrix $\PP^*$ is upper-triangular, thus $\{1-p_1,\ldots,1-p_{M-1},1\}$ are its eigenvalues. Because of the relation 
(\ref{eq:duality}) these are also the eigenvalues of $\PP$.
 \section*{Acknowledgements}
 The author thanks anonymous reviewers for thorough reviews and appreciates the comments 
 and suggestions, which contributed to improving the quality of the publication.
 
%

\bibliographystyle{abbrv}
\bibliography{library}

\begin{thebibliography}{10}

\bibitem{Aldous1986}
D.~Aldous and P.~Diaconis.
\newblock {Shuffling cards and stopping times}.
\newblock {\em American Mathematical Monthly}, 93(5):333--348, 1986.

\bibitem{Aldous1987}
D.~Aldous and P.~Diaconis.
\newblock {Strong Uniform Times and Finite Random Walks}.
\newblock {\em Advances in Applied Mathematics}, 97:69--97, 1987.

\bibitem{Basu2017}
R.~Basu, J.~Hermon, and Y.~Peres.
\newblock {Characterization of cutoff for reversible Markov chains}.
\newblock {\em Annals of Probability}, 45(3):1448--1487, 2017.

\bibitem{Chen2008}
G.~Y. Chen and L.~Saloff-Coste.
\newblock {The cutoff phenomenon for ergodic Markov processes}.
\newblock {\em Electronic Journal of Probability}, 13:26--78, 2008.

\bibitem{Chen2015}
G.-Y. Chen and L.~Saloff-Coste.
\newblock {Computing cutoff times of birth and death chains}.
\newblock {\em Electronic Journal of Probability}, 20:1--47, 2015.

\bibitem{Choi2016}
M.~C.~H. Choi and P.~Patie.
\newblock {A Sufficient Condition for Continuous-Time Finite Skip-Free Markov
  Chains to Have Real Eigenvalues}.
\newblock {\em In: B{\'{e}}lair J., Frigaard I., Kunze H., Makarov R., Melnik
  R., Spiteri R. (eds) Mathematical and Computational Approaches in Advancing
  Modern Science and Engineering.}, pages 529--536, 2016.

\bibitem{Connor2010}
S.~B. Connor.
\newblock {Separation and coupling cutoffs for tuples of independent Markov
  processes}.
\newblock {\em Latin American Journal of Probability and Mathematical
  Statistics}, 7(3):65--77, 2010.

\bibitem{Diaconis1990a}
P.~Diaconis and J.~A. Fill.
\newblock {Strong stationary times via a new form of duality}.
\newblock {\em The Annals of Probability}, 18(4):1483--1522, 1990.

\bibitem{Diaconis2009b}
P.~Diaconis and L.~Miclo.
\newblock {On Times to Quasi-stationarity for Birth and Death Processes}.
\newblock {\em Journal of Theoretical Probability}, 22(3):558--586, jun 2009.

\bibitem{DiaSal_metropolis98}
P.~Diaconis and L.~Saloff-Coste.
\newblock {What do we know about metropolis algorithm?}
\newblock {\em Journal of Computer and System Sciences}, 57:20--36, 1998.

\bibitem{Diaconis2006}
P.~Diaconis and L.~Saloff-Coste.
\newblock {Separation cut-offs for birth and death chains}.
\newblock {\em The Annals of Applied Probability}, 16(4):2098--2122, 2006.

\bibitem{Diaconis1981a}
P.~Diaconis and M.~Shahshahani.
\newblock {Generating a random permutation with random transpositions}.
\newblock {\em Zeitschrift fur Wahrscheinlichkeitstheorie und Verwandte
  Gebiete}, 57(2):159--179, 1981.

\bibitem{Ding2010}
J.~Ding, E.~Lubetzky, and Y.~Peres.
\newblock {Total variation cutoff in birth-and-death chains}.
\newblock {\em Probability Theory and Related Fields}, 146:61--85, 2010.

\bibitem{Doumas2016}
A.~V. Doumas and V.~G. Papanicolaou.
\newblock {The Coupon Collector's Problem Revisited: Generalizing the Double
  Dixie Cup Problem of Newman and Shepp}.
\newblock {\em ESAIM: Probability and Statistics}, 20:367--399, 2016.

\bibitem{Erdos1961}
P.~L. Erdős and A.~R{\'{e}}nyi.
\newblock {On a classical problem of probability theoryd}.
\newblock {\em Publ. Math. Inst. Hung. Acad. Sci.}, Ser. A(6):215--220, 1961.

\bibitem{Feller1971}
W.~Feller.
\newblock {\em {An Introduction to Probability Theory and its Applications,
  Volume 2}}.
\newblock John Wiley {\&} Sons, 2nd edition, 1971.

\bibitem{Fill1996}
J.~A. Fill.
\newblock {An exact formula for the move-to-front rule for self-organizing
  lists}.
\newblock {\em Journal of Theoretical Probability}, 9(1):113--160, 1996.

\bibitem{Fill2009a}
J.~A. Fill.
\newblock {On hitting times and fastest strong stationary times for skip-free
  and more general chains}.
\newblock {\em Journal of Theoretical Probability}, 22(3):587--600, 2009.

\bibitem{Fill2009}
J.~A. Fill.
\newblock {The Passage Time Distribution for a Birth-and-Death Chain: Strong
  Stationary Duality Gives a First Stochastic Proof}.
\newblock {\em Journal of Theoretical Probability}, 22(3):543--557, 2009.

\bibitem{FilLyz14}
J.~A. Fill and V.~Lyzinski.
\newblock {Strong Stationary Duality for diffusion processes}.
\newblock {\em arXiv}, pages 1--33, 2014.

\bibitem{Hermon2016}
J.~Hermon, H.~Lacoin, and Y.~Peres.
\newblock {Total variation and separation cutoffs are not equivalent and
  neither one implies the other}.
\newblock {\em Electronic Journal of Probability}, 21(44):1--36, 2016.

\bibitem{Holst2001}
L.~Holst.
\newblock {Extreme Value Distributions for Random Coupon Collector and Birthday
  Problems}.
\newblock {\em Extremes}, 4(2):129--145, 2001.

\bibitem{Karlin1959}
S.~Karlin and J.~McGregor.
\newblock {Coincidence properties of birth and death processes.}
\newblock {\em Pacific Journal of Mathematics}, 9(4):1109--1140, 1959.

\bibitem{Keilson71b}
J.~Keilson.
\newblock {Log-Concavity and Log-Convexity in Passage Time Densities of
  Diffusion and Birth-Death Processes}.
\newblock {\em Journal of Applied Probability}, 8(2):391--398, 1971.

\bibitem{Lacoin2016}
H.~Lacoin.
\newblock {The Cutoff profile for the Simple-Exclusion process on the circle}.
\newblock {\em Ann. Probab.}, 44(5):3399--3430, 2016.

\bibitem{LevPerWil_mixing_sec_ed}
D.~Levin, Y.~Peres, and E.~Wilmer.
\newblock {\em {Markov Chains and Mixing Times, second edition}}.

\bibitem{2015Lorek_gambler}
P.~Lorek.
\newblock {Generalized Gambler's Ruin Problem: Explicit Formulas via Siegmund
  Duality}.
\newblock {\em Methodology and Computing in Applied Probability},
  19(2):603--613, 2017.

\bibitem{Lorek2016_Siegmund_duality}
P.~Lorek.
\newblock {Siegmund duality for Markov chains on partially ordered state
  spaces}.
\newblock {\em Probability in the Engineering and Informational Sciences},
  pages 1--27, 2017.

\bibitem{LorMar_mon_req}
P.~Lorek and P.~Markowski.
\newblock {Monotonicity requirements for efficient exact sampling with Markov
  chains.}
\newblock {\em Markov Processes And Related Fields}, 23:485--514, 2017.

\bibitem{Lorek2012d}
P.~Lorek and R.~Szekli.
\newblock {Strong stationary duality for M{\"{o}}bius monotone Markov chains}.
\newblock {\em Queueing Systems}, 71(1-2):79--95, mar 2012.

\bibitem{Lubetzky2013}
E.~Lubetzky and A.~Sly.
\newblock {Cutoff for the Ising model on the lattice}.
\newblock {\em Inventiones Mathematicae}, 191(3):719--755, 2013.

\bibitem{MaoZhang2016}
Y.-h. Mao, C.~Zhang, and Y.-h. Zhang.
\newblock {Separation cutoff for upward skip-free chains}.
\newblock {\em Journal of Applied Probability}, 1:299--306, 2016.

\bibitem{Miclo2010a}
L.~Miclo.
\newblock {On absorption times and Dirichlet Eigenvalues}.
\newblock {\em ESAIM: Probability and Statistics}, 14:117--150, 2010.

\bibitem{Neal2008}
P.~Neal.
\newblock {The generalised coupon collector problem}.
\newblock {\em Journal of Applied Probability}, 45(3):621--629, 2008.

\bibitem{Newman60}
D.~Newman.
\newblock {The double dixie cup problem}.
\newblock {\em American Mathematical Monthly}, 67(1):58--61, 1960.

\bibitem{Pak2001}
I.~Pak and V.~H. Vu.
\newblock {On mixing of certain random walks, cutoff phenomenon and sharp
  threshold of random matroid processes}.
\newblock {\em Discrete Applied Mathematics}, 110(2-3):251--272, 2001.

\bibitem{Rota64}
G.-C. Rota.
\newblock {On the foundations of combinatorial theory I. Theory of M{\"{o}}bius
  functions}.
\newblock {\em Probability Theory and Related Fields}, 368:340--368, 1964.

\bibitem{Siegmund1976}
D.~Siegmund.
\newblock {The Equivalence of Absorbing and Reflecting Barrier Problems for
  Stochastically Monotone Markov Processes}.
\newblock {\em The Annals of Probability}, 4(6):914--924, 1976.

\end{thebibliography}

\end{document}